\documentclass[10pt]{amsart}
\usepackage{color}

\usepackage{amsmath,amssymb}

\definecolor{c20}{rgb}{0.,0.7,0.}
\definecolor{c30}{rgb}{0.,0.,1.}
\definecolor{c40}{rgb}{1,0.1,0.7}
\definecolor{c50}{rgb}{1,0,0}
\definecolor{c60}{rgb}{1,0.9,0.1}

\newcommand{\abs}[1]{\left\lvert #1 \right\rvert}

\newcommand{\norm}[1]{\lVert #1 \rVert}

\newcommand{\E}[1]{\mathbb{E}\left\{#1\right\}}

\newcommand{\pk}[1]{\mathbb{P} \left\{ #1 \right \} }

\newcommand{\R}{\mathbb{R}}

\newcommand{\BQN}{\begin{eqnarray}}
\newcommand{\EQN}{\end{eqnarray}}
\newcommand{\BQNY}{\begin{eqnarray*}}
\newcommand{\EQNY}{\end{eqnarray*}}

\newcommand{\BS}{\begin{sat}}
\newcommand{\ES}{\end{sat}}
\newcommand{\BT}{\begin{theo}}
\newcommand{\ET}{\end{theo}}
\newcommand{\BK}{\begin{korr}}
\newcommand{\EK}{\end{korr}}

\newcommand{\BD}{\begin{de}}
\newcommand{\ED}{\end{de}}

\newcommand{\BIT}{\begin{itemize}}
\newcommand{\EIT}{\end{itemize}}
\newcommand{\BDI}{\begin{description}}
\newcommand{\EDI}{\end{description}}

\newcommand{\BRM}{\begin{remarks}}
\newcommand{\ERM}{\end{remarks}}

\newcommand{\BEL}{\begin{lem}}
\newcommand{\EEL}{\end{lem}}

\newtheorem{theo}{Theorem}[section]
\newtheorem{sat}[theo]{Proposition}
\newtheorem{de}[theo]{Definition}
\newtheorem{lem}[theo]{Lemma}

\newtheorem{example}[theo]{Example}
\newtheorem{korr}[theo]{Corollary}
\newtheorem{remark}[theo]{Remark}
\newtheorem{remarks}[theo]{Remarks}

\newcommand{\nelem}[1]{{Lemma \ref{#1}}}

\newcommand{\netheo}[1]{{Theorem \ref{#1}}}

\newcommand{\prooftheo}[1]{ \textsc{\bf Proof of Theorem} \ref{#1}:}

\newcommand{\prooflem}[1]{\textsc{\bf Proof of Lemma} \ref{#1}:}

\newcommand{\COM}[1]{}

\newcommand{\QED}{\hfill $\Box$}

\newcommand{\kb}[1]{\boldsymbol{#1}}
\newcommand{\vk}[1]{\kb{#1}}

\def\vn{\varepsilon}

\def\IF{\infty}

\def\LT{\left}
\def\RT{\right}

\def\HH{\mathcal{H}}

\def\vn{\varepsilon}

\def\Var{\mathrm{Var}}

\def\dsig{\dot{\sigma}}

\def\mlG{\mathcal{G}}

\def\to{\rightarrow}

\def\wtB{\widetilde B}
\def\wtH{\widetilde H}
\def\wtT{\widetilde {\vk T}}

\def\mcB{\mathcal B}
\def\barR{\overline{\R}}
\def\vto{\overset{v}\to}

\def\isig1u{\overleftarrow{\sigma_1}(u)}

\def\MT{\mathcal{T}}
\def\1Hi{\frac{1}{H_i}}

\def\cL#1{\textcolor{red}{#1}}
\def\cL#1{\textcolor{black}{#1}}

\begin{document}

\title{Extrema of multi-dimensional Gaussian processes over random intervals }

 \author{Lanpeng Ji}
\address{Lanpeng Ji, School of Mathematics, University of Leeds, Woodhouse Lane, Leeds LS2 9JT, United Kingdom
}
\email{l.ji@leeds.ac.uk}

\author{Xiaofan Peng}
\address{Xiaofan Peng, School of Mathematical Sciences, University of Electronic Science and Technology of China, Chengdu 610054, China}
\email{xfpeng@uestc.edu.cn}

\bigskip

\date{\today}
 \maketitle

 {\bf Abstract:} This paper studies the joint tail asymptotics of extrema of the multi-dimensional Gaussian process  over random intervals defined as
 $$
P(u):=\mathbb{P}\left\{\cap_{i=1}^n \left(\sup_{t\in[0,\mathcal{T}_i]} ( X_{i}(t) +c_i t )>a_i u \right)\right\}, \ \ \ u\to\infty,
$$
 where $X_i(t), t\ge0$, $i=1,2,\cdots,n,$ are independent centered Gaussian processes with stationary increments, $\boldsymbol{\mathcal{T}}=(\mathcal{T}_1, \cdots, \mathcal{T}_n)$ is a regularly varying random vector with positive components, which is independent of  the Gaussian processes, and $c_i\in \mathbb{R}$, $a_i>0$, $i=1,2,\cdots,n$. Our result shows that the structure of the asymptotics of $P(u)$ is determined  by the signs of the drifts $c_i$'s.
We also discuss a relevant multi-dimensional regenerative model and derive the corresponding ruin probability.
\medskip

%Generalized Pickands constant.\\


\section{Introduction }

Let $X(t), t\ge0$  be an almost surely (a.s.) continuous centered Gaussian process with stationary increments and $X(0)=0$. Motivated by its applications to the hybrid fluid and ruin models, the seminal paper \cite{debicki2004supremum} derived   the exact tail asymptotics of
\BQN\label{eq:XXt}
\pk{\sup_{t\in [0,\mathcal T]} X(t) >u},\ \ \ u\to\IF,
\EQN
with $\mathcal T$ being an independent of $X$ regularly varying random variable. Since then the study of the tail asymptotics of supremum on random interval has attracted substantial interest in the literature. %;  see the seminal papers \cite{ZBD05, debicki2004supremum}.
We refer to \cite{arendarczyk2011asymptotics, arendarczyk2011exact, TH13, DEJ14, Are17, DP20} for various extensions to general (non-centered) Gaussian  or Gaussian-related processes. In the aforementioned contributions, various different tail distributions for $\mathcal T$ have been discussed, % (e.g., regularly varying tail distribution or  Weibullian-type tail distribution),
and it has been shown that the  variability of $\mathcal T$  influences the form of the asymptotics of \eqref{eq:XXt}, leading to qualitatively different structures.

The primary aim of this paper is to analyze the asymptotics of a multi-dimensional counterpart of \eqref{eq:XXt}. More precisely, consider a multi-dimensional centered Gaussian process
\BQN \label{eq:vkX}
\vk X(t)=(X_1(t), X_2(t), \cdots, X_n(t)),\ \ \ t\ge0,
\EQN
with independent coordinates,  each $X_i(t), t\ge0,$ has stationary increments, a.s. continuous sample paths and $X_i(0)=0$, and
% the standard deviation functions $\sigma_i(t)=\sqrt{Var(X_i(t))}$ satisfies the assumptions {\bf S1}-{\bf S5} with the parameters involved indexed by $i$; e.g., {\bf S1} is satisfied with $H_i\in (0,1)$.
let  $ \vk{\mathcal{ T}}=(\mathcal{T}_1, \cdots, \mathcal{T}_n)$ be a regularly varying random vector with positive components, % with index $\alpha$ and limiting measure $\nu$
which is independent of $\vk X$.
We are interested in the exact asymptotics of
\BQN\label{eq:Pu}
P(u):=\pk{\cap_{i=1}^n \LT(\sup_{t\in[0,\mathcal{T}_i]} ( X_{i}(t) +c_i t )>a_i u \RT)}, \ \ \ u\to\IF,
\EQN
where $c_i\in \R$, $a_i>0$, $i=1,2,\cdots,n$.

Extremal analysis of multi-dimensional Gaussian processes has been an active research area in recent years; see \cite{Debicki10, DHJT15, AP19, Pham20, DHW20} and references therein. In most of these contributions, the  asymptotic behaviour of the probability that $\vk X$ (possibly with trend) enters an upper orthant over a finite-time or infinite-time interval is discussed, this problem is also connected with the conjunction problem for Gaussian processes firstly studied by Worsley and Friston \cite{WF00}. Investigations on the joint tail asymptotics of multiple extrema as defined in \eqref{eq:Pu} have been known to be more challenging. Current literature has only focused on the case  with deterministic times $\mathcal{T}_1=\cdots=\mathcal{T}_n$ and some additional assumptions on the correlation structure of $X_i$'s. In  \cite{MR2145669,  Debicki10} large deviation type results are obtained, and more recently in  \cite{DJT20, DHK20} exact asymptotics are obtained for correlated two-dimensional Brownian motion. It is worth mentioning that a large derivation  result for the multivariate maxima of  a discrete Gaussian model is discussed recently in \cite{HH19}.

In order to avoid more technical difficulties,  the coordinates of the multi-dimensional process $\vk X$ in   \eqref{eq:vkX} are assumed to be independent. The dependence among the extrema in \eqref{eq:Pu} is driven by the structure of the multivariate regularly varying $\vk {\mathcal{T}}$.  Interestingly, we observe in Theorem \ref{Thm:Main} that the form of the asymptotics of \eqref{eq:Pu} is determined by the signs of the drifts $c_i$'s.

Apart from its theoretical interest, the motivation to analyse the asymptotic properties of $P(u)$ is related to numerous applications in modern multi-dimensional risk theory, financial mathematics or fluid queueing networks.
For example, we consider an insurance company which runs $n$  lines of business. The surplus process of the $i$th business line can be  modelled by a time-changed Gaussian process
\BQNY
R_i(t)=a_i u+c_i Y_i(t) - X_i(Y_i(t)),\ \ t\ge0, % i=1,\cdots, n.
\EQNY
where $a_i u>0$ is the initial capital (considered as a proportion of $u$ allocated to the $i$th business line, with $\sum_{i=1}^n a_i=1$), $c_i>0$ is the net premium rate, $X_i(t), t\ge0$ is the net loss process, and $Y_i(t), t\ge 0$ is a positive increasing function modelling the so-called  ``operational time" for the $i$th business line. We refer to \cite{AsmAlb10, JR18} and \cite{DEJ14} for detailed discussions on multi-dimensional risk models and time-changed  risk models, respectively. Of interest in risk theory is the study of the  probability of ruin of all the business lines   within some finite (deterministic) time $T>0$, defined by
\BQNY
\varphi(u):=\pk{\cap_{i=1}^n \LT(\inf_{t\in[0,T]} R_i(t)<0 \RT)} =\pk{\cap_{i=1}^n \LT(\sup_{t\in[0,T]} ( X_{i}(Y_i(t)) +c_i Y_i(t) )>a_i u \RT)}.
\EQNY
If additionally all the operational time  processes $Y_i(t),t\ge0$ have a.s. continuous sample paths, then we have $\varphi(u) = P(u)$ with $\vk {\mathcal{T}}=\vk{Y}(T)$, and thus the derived result %in Theorem \ref{Thm:Main}
can be applied to estimate this ruin probability. Note that the dependence among different business lines is introduced by the dependence among the operational time  processes $Y_i$'s. %extrema in \eqref{eq:Pu} is mainly driven by the structure of the multivariate regularly varying $\vk {\mathcal{T}}$.
As a simple example we can consider $Y_i(t)= \Theta_i t,$ $ t\ge 0$, with $\vk \Theta= (\Theta_1,\cdots, \Theta_n)$ being a multivariate regularly varying random vector.
%Note that the dependence among the extrema in \eqref{eq:Pu} is mainly driven by the structure of the multivariate regularly varying $\vk {\mathcal{T}}$.
%Normally, these business time processes $Y_i, i=1,\cdots, n,$ are correlated.
Additionally,  multi-dimensional time-changed (or  subordinate)  Gaussian processes have been  recently proved to be good candidates for modelling the log-return processes of multiple assets; see, e.g., \cite{BPS01, LS10, Kim12}. As    the joint distribution of extrema of asset returns is  important in finance problems, e.g., \cite{HKR98}, we expect the obtained results for \eqref{eq:Pu} might also be interesting in financial mathematics.

As a relevent application, we shall discuss a multi-dimensional regenerative model, which is motivated from its relevance to risk model and fluid queueing model. Essentially, the multi-dimensional regenerative process is a process with a random alternating environment, where an independent multi-dimensional fractional Brownian motion (fBm) with trend is assigned at each  environment alternating time. We refer to Section 4 for more detail. By analysing a related multi-dimensional perturbed random walk,  we obtain in Theorem \ref{Thm:Qu} the ruin probability of the multi-dimensional regenerative model. This generalizes some of the results in \cite{PZ07} and \cite{ZBD05}  to the multi-dimensional setting. Note in passing that some related stochastic models with random sampling or resetting have been discussed in the recent literature; see, e.g., \cite{KellaWhitt, Corina20, Rat20}.

\medskip

{\it Organization of the rest of the paper}: In Section 2 we introduce some notation, recall the definition of the multivariate  regularly variation, and present some preliminary results on the extremes of
one-dimensional Gaussian  process. The   result  for \eqref{eq:Pu} is displayed in Section 3, and the ruin probability of the multi-dimensional regenerative model is discussed in Section 4. The proofs  are relegated to  Section 5 and Section 6.  Some useful results on multivariate regularly variation are discussed in the Appendix.

%%%%%%%%%%%%%%%%%%%%%%%%%%%%%%%%%

\section{Notation and Preliminaries}

We shall use some standard notation which is common when dealing with vectors.
All the operations on vectors are meant componentwise. For instance, for any given
 $ \vk{x} = (x_1,\ldots,x_n)\in \R ^n$ and $\vk{y} = (y_1,\ldots,y_n) \in \R ^n $, we write $\vk{x} \vk y = (x_1y_1, \cdots, x_ny_n)$,  and write $ \vk{x} > \vk{y} $ if and only if
  $ x_i > y_i $  for all $ 1 \leq i \leq n $. %,  and write $\lambda\vk{y}= (\lambda y_1 \ldot \lambda y_n)$ for any $\lambda\in \R, \vk y\in \R^n$.
%Further we set  $ \vk{0}  := (0,\ldots,0)\ccP{\in\R^n} $ and $ \vk{1} : = (1,\ldots,1)\ccP{\in\R^n}$. \\
Furthermore, for two positive functions $f, h$ and some $u_0\in[-\IF, \IF]$, write $ f(u) \lesssim h(u)$ or $h(u)\gtrsim  f(u)$ if $ \limsup_{u \to u_0} f(u) /h(u)  \le 1 $, write $h(u)\sim  f(u)$ if $ \lim_{u \to u_0} f(u) /h(u)  = 1 $,  write $ f(u) = o(h(u)) $ if $ \lim_{u \to u_0} {f(u)}/{h(u)} = 0$, and write $ f(u) \asymp h(u) $ if $ f(u)/h(u)$ is bounded from both below and above for all sufficiently large $u$.

Next, let us recall the definition and some implications of multivariate regularly variation. We  refer to \cite{Jessen06, Res07, Hult06} for more detailed discussions. Let  $\barR_0^n=\barR^n \setminus \{\vk 0\}$ with $\barR=\R\cup\{-\IF, \IF\}$.
%{\bf Definition: }
An $\R^n$-valued random vector $\vk X$ is said to be {\it regularly varying} if there exists a non-null Radon measure $\nu$ on the Borel $\sigma$-field $\mcB(\barR_0^n)$  with  $\nu(\barR^n  \setminus \R^n)=0$  such that
\BQNY
\frac{\pk{x^{-1} \vk X\in \cdot} }{\pk{\abs{\vk X} >x}}\ \vto \ \nu(\cdot),\ \ \ \ x\to \IF.
\EQNY
Here $\abs{\ \cdot\ }$ is any norm in $\R^n$ and $\vto$ refers to vague convergence on $\mcB(\barR_0^n)$.
 It is known that $\nu$ necessarily satisfies the homogeneous property
$\nu(s K) =s^{-\alpha} \nu (K)$, $s>0$, for some $\alpha>0$ and all Borel set $K$ in $\mcB(\barR_0^n)$.
In what follows, we say that such defined $\vk X$ is regularly varying with index $\alpha$ and limiting measure $\nu$.
%If the limit measure is irrelevant we also write $\vk X \in MRV(\alpha)$.
An implication of the homogeneity property of $\nu$ is that all the rectangle sets of the form $[\vk a, \vk b]=\{\vk x: \vk a \le \vk x\le \vk b\}$ in $\barR_0^n$ are $\nu$-continuity sets. %, where $\vk a<\vk b$ with $<, \le$ defined in the natural componentwise way.
Furthermore, we have that $\abs{\vk X}$ is regularly varying at infinity with index $\alpha$, i.e., $\pk{\abs{\vk X}>x}\sim x^{-\alpha} L(x), x\to\IF$, with some slowly varying function $L(x)$. Some useful results on the multivariate regularly variation %varying $\vk{\mathcal{T}}$
are discussed in Appendix.

%Throughout the paper, all the vectors with bold letters are considered to be $n$-dimensional.

In what follows, we review some results on the extremes of one-dimensional Gaussian  process with nagetive drift derived in \cite{dieker2005extremes}.
Let $X(t), t\ge0 $  be an a.s. continuous centered Gaussian process with stationary increments and $X(0)=0$, and let $c>0$ be some constant. We shall present  the exact asymptotics of
$$\psi (u):=\pk{\sup_{t\ge0} (X(t)-ct)>u},\ \ \ u\to\IF.$$

Below are some assumptions that the variance function $\sigma^2(t)=\Var(X(t))$ might satisfy:
\begin{itemize}

\item[\bf C1:] $\sigma$ is continuous on $[0,\IF)$ and ultimately strictly increasing;

\item[\bf C2:] $\sigma$ is regularly varying at infinity with index $H$ for some $H\in(0,1)$;

\item[\bf C3:] $\sigma$ is regularly varying at 0 with index $\lambda$ for some $\lambda\in(0,1)$;

\item[\bf C4:] %$\sigma^2$ is ultimately continuously differentiable and its first derivative $\dsig^2$ is ultimately monotone;
$\sigma^2$ is ultimately twice continuously differentiable and its first derivative $\dsig^2$ and second derivative $\ddot{ \sigma}^2$ are both ultimately monotone.

\end{itemize}
%There exist some $G ,\gamma>0$ such that
%\BQNY
%\sigma^2(t)\le G t^\gamma
%\EQNY
% holds on a neighborhood of zero;

Note that in the above  $\dsig^2$ and $\ddot{ \sigma}^2$ denote the first and second derivative of $\sigma^2$, not the square of the derivatives of $\sigma$.
 In the sequel, provided it exists we denote by $\overleftarrow{\sigma}$ an asymptotic inverse near infinity or zero of $\sigma$; recall that it is (asymptotically uniquely) defined by $\overleftarrow{\sigma}(\sigma(t))\sim \sigma(\overleftarrow{\sigma}(t))\sim t.$  It depends on the context whether  $\overleftarrow{\sigma}$ is an asymptotic inverse near zero or infinity. % (need to check if we used the monotonicity of $\sigma$ in the proof)

One known example that satisfies the assumptions {\bf C1}-{\bf C4}  is the  fBm  $\{B_H(t), t\ge 0\}$  with Hurst index $H\in(0,1)$, i.e., an $H$-self-similar centered Gaussian process with stationary increments and covariance function given by
\BQNY
Cov(B_H(t),B_H(s))=\frac{1}{2}(\abs{t}^{2H}+\abs{s}^{2H}-\abs{t-s}^{2H}),\ \ t,s\in \R.
\EQNY

We introduce the following notation:
\BQNY
C_{H,\lambda_1,\lambda_2}=\sqrt{2^{1-1/\lambda_2}\pi} \lambda_1 \LT(\frac{1}{H}\RT)^{1/\lambda_2}\LT(\frac{H}{1-H}\RT)^{\lambda_1+ H -\frac{1}{2}+\frac{1}{\lambda_2}\LT(1- H \RT)}.
\EQNY
For  an a.s. continuous centered Gaussian process $Z(t), t\ge 0$ with stationary increments and variance function $\sigma_Z^2$, we define the generalized Pickands constant
\BQNY
\HH_Z=\lim_{T\to\IF}\frac{1}{T}\E{\exp\LT(\sup_{t\in[0,T]}(\sqrt{2}Z(t)-\sigma_Z^2(t))\RT)}
\EQNY
provided both the expectation and the limit exist. When $Z=B_H$ the constant $\HH_{B_H}$ is the well-known Pickands constant; see \cite{Pit96}. %; see [] for recent studies on Pickands constant and related generalizations.
For convenience, sometimes we also write $\HH_{\sigma_Z^2}$ for $\HH_Z$. Denote in the following by $\Psi(\cdot)$ the survival function of the $N(0,1)$ distribution. It is known that
\BQN \label{eq:Psiu}
\Psi(u)=\frac{1}{\sqrt{2\pi}} \int_u^\IF e^{-\frac{x^2}{2}} dx\  \sim \ \frac{1}{\sqrt{2\pi} u}  e^{-\frac{u^2}{2}},\ \ \ \ u\to\IF.
\EQN

 The following result is derived in Proposition 2 in \cite{dieker2005extremes} (here we consider a particular trend function  $\phi(t)= ct, t\ge 0$). %For simplicity we assume $A=1$ (will be removed  later on).
\BS\label{CorXc}
Let $X(t),  t\ge0$  be an a.s. continuous centered Gaussian process with stationary increments and $X(0)=0$. Suppose that {\bf C1--C4} hold. We have, as $u\to\IF$

(i) if $\sigma^2(u)/u\to \IF$, then
\BQNY
\psi(u) \sim \HH_{B_H} C_{H,1,H} \LT(\frac{1-H}{H}\RT)\frac{c^{1-H}\sigma(u) }{\cL{\overleftarrow{\sigma}}(\sigma^2(u)/u)}\Psi\LT(\inf_{t\ge0}\frac{u(1+t)}{\sigma(ut/c)}\RT);
\EQNY

(ii) if $\sigma^2(u)/u\to \mlG \in(0,\IF)$, then
\BQNY
\psi(u) \sim \HH_{(2c^2/\mlG^2 )\sigma^2}  \LT(\frac{\sqrt{2/\pi}}{c^{1+H}H}\RT) \sigma(u)\Psi\LT(\inf_{t\ge0}\frac{u(1+t)}{\sigma(ut/c)}\RT);
\EQNY

(iii) if $\sigma^2(u)/u\to 0$, then [here we need regularity of $\sigma$ and its inverse at 0]
\BQNY
\psi(u) \sim \HH_{B_\lambda} C_{H,1,\lambda} \LT(\frac{1-H}{H}\RT)^{H/\lambda}\frac{c^{\cL{-1-H+2H/\lambda}}\sigma(u) }{\cL{\overleftarrow{\sigma}}(\sigma^2(u)/u)}\Psi\LT(\inf_{t\ge0}\frac{u(1+t)}{\sigma(ut/c)}\RT).
\EQNY
\ES

\COM{
\begin{remark}
Note that the result presented in \cite{dieker2005extremes} is for more general trend function while here the trend function is $\phi(t)= ct, t\ge 0$. Furthermore, condition {\bf C3} is not needed for (i)-(ii) to hold, and for (iii) this can be replaced by the following assumption: There exist some $G ,\gamma>0$ such that
$\sigma^2(t)\le G t^\gamma $ holds on a neighbourhood of zero. For simplicity, we keep the slightly stricter condition  {\bf C3}.
\end{remark}
}

\COM{ %%%%%%%
Note that in the proof of the main results in \cite{dieker2005extremes} the minimum point $t_u^*$ of the function
\BQNY
f_u(t) :=\frac{u(1+t)}{\sigma(ut/c)}, \ \ t\ge 0,
\EQNY
plays an important role. It has been discussed therein that $t_u^*$ converges, as $u\to\IF,$ to $t^*:=H/(1-H)$ which is the unique minimum point of $\lim_{u\to\IF} f_u(t) \sigma(u)/u= (1+t)/ (t/c)^H, t\ge 0$. In this sense, $t_u^*$ is asymptotically unique. Moreover, we have the following implications, the proof of it is displayed in Section \ref{Sec:pro}.
\BS \label{Prop:2}
For any fixed $0<\vn<t^*/c$, we have, as $u\to\IF,$
\BQNY
\pk{\sup_{t\in 0, (t^*/c +\vn) u]} ( X(t)-ct) >u}\sim P_c(u). % ,\ \ \ \pk{\sup_{t\in[0, (t^*/c -\vn) u]} ( X(t)-ct) >u} =o(P_c(u)).
\EQNY
Furthermore, we have that for any $\gamma>0$
\BQNY
\lim_{u\to\IF}\frac{\pk{  \sup_{t\in[0,(t^*/c-\vn) u ]} X(t) -c t >u    } }{P_c(u) u^{-\gamma}}=0.
\EQNY

\ES

} %%%%%%%5

As a special case of the Proposition \ref{CorXc} we have the following result (see Corollary 1 in \cite{dieker2005extremes} or \cite{JR18}). This will be useful in the proofs below.

\BK\label{eq:psi_i}
If $X(t)=B_H(t), t\ge 0$ is the fBm with index $H\in(0,1)$, then  as $u\to\IF$
\BQNY
 \pk{\sup_{t\ge 0} (B_{H}(t)-c t)>  u} %\nonumber\\
 \sim K_H  \mathcal{H}_{{B_H}}    u^{H+1/H-2}\ \! \Psi\left(\frac{c^{H}  u^{1-H}}{H^{H} (1-H)^{1-H}}\right).
\EQNY
with constant $K_H=2^{\frac{1}{2}-\frac{1}{2H}}\frac{\sqrt{\pi}}{\sqrt{H(1-H)}}   \left(\frac{c^{H}  }{H^{H} (1-H)^{1-H}}\right)^{1/H-1}$.
%and
%$\mathcal{T}_i^*=\frac{H_i}{(H_i-1) c_i}.$
\EK

\section{Main results }

Without loss of generality, we assume that \cL{in \eqref{eq:Pu} there are $n_-$ coordinates with negative drift, $n_0$ coordinates without drift and $n_+$ coordinates with positive drift, i.e.,}
\BQNY
&&c_i<0, \ \ i=1,\cdots, n_-, \\
&&c_i=0, \ \ i=n_-+1,\cdots, n_-+n_0,\\
&&c_i>0, \ \ i=n_-+n_0+1,\cdots,  n,
\EQNY
where \cL{$0\le n_-,n_0,n_+\le n$ such that $n_-+n_0+n_+=n$.}
We impose the following assumptions for the standard deviation functions $\sigma_i(t)=\sqrt{Var(X_i(t))}$ of the Gaussian processes $X_i(t), i=1,\cdots,n$.

{\bf Assumption I:} For $i=1,\cdots, n_-,$  $\sigma_i(t)$ satisfies the assumptions {\bf C1}-{\bf C4} with the parameters involved indexed by $i$.
 For $i=n_-+1,\cdots, n_-+n_0,$  $\sigma_i(t)$ satisfies the assumptions {\bf C1}-{\bf C3} with the parameters involved indexed by $i$.  For $i=n_-+n_0+1,\cdots, n,$  $\sigma_i(t)$ satisfies the assumptions {\bf C1}-{\bf C2} with the parameters involved indexed by $i$.

Denote
\BQN\label{eq:xieta}
\xi_i:=\sup_{t\in[0,1]} B_{H_i}(t),\ \ \ %\eta_i:=B_{H_i}(1)
t^*_i=\frac{H_i}{1-H_i}.
\EQN
Given a Radon measure $\nu$, define
\BQN\label{def:nutil}
\widetilde{\nu}(K)=:\E{\nu(\vk \xi^{-1/{\vk H}} K)}, \quad K\subset  \mathcal{B}([0,\IF]^n \setminus \{\vk 0\}),
 \EQN
where $\vk \xi^{-1/{\vk H}} K=\{(\xi_1^{-1/{H_1}}d_1,\cdots, \xi_n^{-1/{H_n}}d_n), (d_1,\cdots,d_n)\in K\}$.
Further, note that for $i=1,\cdots, n_-,$ (where $c_i<0$), the asymptotic  formula, as $u\to\IF$, of
\BQN\label{def:psiiu}
\psi_i(u)=\pk{\sup_{t\ge 0} (X_i(t)+c_i t)>u}.
\EQN
  is available from  Proposition \ref{CorXc} under Assumption {\bf I}.

\BT \label{Thm:Main} Suppose that  $\vk X(t), t\ge 0$ satisfies the Assumption  \textbf{I}, and   $\vk{\mathcal{T}}$ is an \cL{ independent of $\vk X$} regularly varying random vector with index $\alpha$ and limiting measure $\nu$. Further assume, without loss of generality, that there are $m(\leq n_0)$ positive constants $k_i$'s such that $\overleftarrow{\sigma_i}(u) \sim k_i \overleftarrow{\sigma}_{n_-+1}(u) $ for $i=n_-+1,\cdots, n_-+m$ and $\overleftarrow{\sigma_i}(u) =o(\cL{\overleftarrow{\sigma}_{n_-+1}}(u))$ for $i=n_-+m+1,\cdots, n_-+n_0$.  We have, with the convention $\prod_{i=1}^{0}=1$,
\begin{itemize}
\item[(i)] If $ n_0>0$, then,  as $u\to\IF$,
\BQNY
P(u)\sim \widetilde{\nu}(( \vk k\vk a_0 ^{1/{H_{n_-+1}}},\vk\IF]) \ \pk{\abs{\vk {\mathcal T}} > \overleftarrow{\sigma}_{n_-+1}(u)}\  \prod_{i=1}^{n_-} \psi_i(a_i u),  %\pk{\cap_{i=n_-+1}^{n_-+n_0} \LT( \mathcal{T}_i \xi_i^{1/H_i}  >a_i^{1/H_i} \cL{\overleftarrow{\sigma_i}}(u)  \RT)}.  %\overleftarrow{\sigma_i}
\EQNY
 where $\widetilde{\nu}$ and $\psi_i's$ are defined in \eqref{def:nutil} and \eqref{def:psiiu}, respectively,
%$\widetilde{\nu}(K)=\E{\nu(\vk \xi^{-1/{\vk H}} K)}$, with $\vk \xi^{-1/{\vk H}} K=\{(\xi_1^{-1/{H_1}}d_1,\cdots, \xi_n^{-1/{H_n}}d_n), (d_1,\cdots,d_n)\in K\}$ for any $K\subset  [0,\IF]^n \setminus \{\vk 0\}$,
and \\
$\vk k\vk a_0^{1/{H_{n_-+1}}}=(0,\cdots, 0, k_{n_-+1}a_{n_-+1}^{1/H_{n_-+1}}, \cdots, k_{n_-+m}a_{n_-+m}^{1/H_{n_-+1}},  0,\cdots, 0).$
%where $D_H $%\{k: \ n_-+1\le k \le n_-+n_0,  H_k =\min_{n_-+1\le i\le n_-+n_0}H_i \}.
%is a set of indexes in $\{1,2,\cdots,n\}$ such that $\cL{\overleftarrow{\sigma_i}}(u) \sim\cL{\overleftarrow{\sigma_j}}(u)$ for any $i,j\in D_H$ and $\cL{\overleftarrow{\sigma_i}}(u) =o(\cL{\overleftarrow{\sigma_j}}(u))$ for any $j\in D_H$, $i\in \{1,2,\cdots,n\}\setminus D_H$.
\item[(ii)]  If $n_0=0$, then,  as $u\to\IF$,
\BQNY
P(u)\sim \nu((\vk a_1, \vk\IF])\ \pk{\abs{\vk {\mathcal T}}>u}\ \prod_{i=1}^{n_-} \psi_i(a_i u),
% \pk{\cap_{i=1}^{n_-} \LT( \mathcal{T}_i> t^*_i/\abs{c_i}  u \RT), \cap_{i=n_-+1}^n \LT( c_i \mathcal{T}_i>a_i u \RT)}.
\EQNY
where $\vk a_1 =( t^*_1/\abs{c_1}\cdots,    t^*_{n_-}/\abs{c_{n_-}} , a_{n_-+1}/{ c_{n_-+1}}, \cdots,  a_n/c_n).$
\end{itemize}

\ET

\begin{remark}
As a special case, we can obtain from Theorem \ref{Thm:Main} some results for the one-dimensional model.
Specifically, let $c>0$ be some constant, then as $u\to\IF$,
 \BQN \label{onedce0}
\pk{ \sup_{t\in[0,\mathcal{T}]}  X(t)> u } &\sim & \E{\LT(\sup_{t\in[0,1]}B_H(t)\RT)^{\alpha/H}}\pk{\mathcal{T}>\overleftarrow{\sigma}(u)},\\
 \label{onedcs0}
\pk{ \sup_{t\in[0,\mathcal{T}]} ( X(t) -c t )> u } &\sim & (c(1-H)/H)^\alpha \pk{\mathcal{T}>u}\psi(u),\\ \label{onedcg0}
\pk{ \sup_{t\in[0,\mathcal{T}]} ( X(t) +c t )> u } &\sim & c^\alpha \pk{\mathcal{T}>u}.
\EQN
Note that \eqref{onedce0} is derived  in Theorem 2.1 of \cite{debicki2004supremum}, \eqref{onedcs0} is discussed in \cite{DEJ14} only for the  fBm case. The result in \eqref{onedcg0} seems to be new.
\end{remark}

We conclude this section with an interesting example of multi-dimensional subordinate Brownian motion; see, e.g., \cite{ LS10}.
\begin{example}\label{Exm}
For each $i=0,1,\cdots, n$, let $\{S_i(t), t\ge0\}$ be independent $\alpha_i$-stable subordinator with $\alpha_i\in(0,1)$,
i.e.,  $S_i(t)\overset{D}=\mathcal S_{\alpha_i}(t^{1/\alpha_i}, 1,0)$, where
$\mathcal S_\alpha(\sigma, \beta, d)$ denotes a stable random variable with stability index  $\alpha$,
scale parameter $\sigma$, skewness parameter $\beta$ and drift parameter $d$. It is known  (e.g., Property 1.2.15 in \cite{ST94}) that for any fixed constant $T>0$,
\BQNY
 \pk{S_i(T)>t}\ \sim\ C_{\alpha_i, T}  t^{-\alpha_i},\quad   \quad  t \rightarrow \infty,
\EQNY
with $
C_{\alpha_i,T}=\frac{T }{\Gamma(1-\alpha_i)\cos(\pi \alpha_i/2)}.
$ Assume $\alpha_0<\alpha_i,$ for all $i=1,2\cdots, n.$ Define an  $n$-dimensional subordinator as
$$
\vk Y(t) :=(S_0(t)+S_1(t), \cdots, S_0(t)+S_n(t)), \ \ \ t\ge 0.
$$
We consider an $n$-dimensional subordinate Brownian motion with drift  defined as
\BQNY
\vk   X(t)= (B_1(Y_1(t))+c_1 Y_1(t), \cdots, B_n(Y_n(t))+c_n Y_n(t)),\ \ \ \ t\ge0,
\EQNY
where $ B_i(t), t\ge0$, $i=1,\cdots, n,$ are independent standard  Brownian motions which are independent of $\vk Y$ and $c_i\in \R$. %Assume for simplicity that $\abs{c_i}$ takes value either 1 or 0.
Define, for any $a_i>0, i=1,2,\cdots, n$,  $T>0$ and $u>0,$
\BQNY
P_B(u):=\pk{\cap_{i=1}^n \LT(\sup_{t\in[0, T]} ( B_{i}(Y_i(t)) +c_i Y_i(t) )> a_i u \RT)}.
\EQNY
For illustrative purpose and to avoid further technicality, we only consider the case where all $c_i$'s in the above have the same sign. As an application of Theorem \ref{Thm:Main} we  obtain  the asymptotic behaviour of $P_B(u), u\to\IF,$ as follows: % (choose $\abs{\ \cdot\ }$ to be the $L^1$ norm):
\begin{itemize}
\item[(i)] If  $c_i>0$ for all $i=1,\cdots, n,$ then
$
  P_B(u) \sim  %\sim  \pk{ S_0(T) > \max_{i=1}^n (a_i/c_i)  u   }  \sim\
  C_{\alpha_0,T} (\max_{i=1}^n (a_i/c_i) u) ^{-\alpha_0}  .
$

\item[(ii)]   If  $c_i=0$ for all $i=1,\cdots, n,$ then
%\BQNY
%\pk{\cap_{i =1}^{n_0} \LT( S(T) (B_i(1))^{2}  >a_i^{2} u^2 \RT)}
%\lesssim P_B(u) \lesssim \pk{\cap_{i =1}^{n_0} \LT(  S(T) \xi_i^{2}  >a_i^{2} u^2 \RT)}.
%\EQNY
%Or better
$ P_B(u)  \asymp u^{-2\alpha_0}.$
\item[(iii)]  If  $c_i<0$ and the density function of $S_i(T)$ is ultimately monotone  for all $i=0, 1,\cdots, n,$   then
$
\ln P_B(u) \sim   2 \sum_{i=1}^n (a_i c_i) u.
$
\end{itemize}

The proof of the above is displayed in Section \ref{Sec:promain}.

\end{example}

\section{Ruin probability of a  multi-dimensional regenerative model}

%Check the multi-dimensional model without noise ...
As it is known in the literature that the maximum of random processes over random interval is relevant to the regenerated models (e.g., \cite{PZ07, ZBD05}),
this section is focused on  a multi-dimensional regenerative model which is motivated from its applications in queueing theory and ruin theory.  More precisely, there are four elements in this model:  Two sequences of strictly positive random variables, $\{T_i: i\ge 1\}$ and  $\{S_i: i\ge 1\}$, and two sequences of $n$-dimensional   processes, $\{\{\vk X^{(i)}(t), t\ge 0\}: i\ge 1\}$ and $\{\{\vk Y^{(i)}(t), t\ge 0\}: i\ge 1\}$, where $\vk X^{(i)}(t)=(X_1^{(i)}(t),\cdots, X_n^{(i)}(t))$ and $\vk Y^{(i)}(t)=(Y_1^{(i)}(t),\cdots, Y_n^{(i)}(t))$. We assume that the above four elements are mutually independent. Here $T_i, S_i$ are two successive times representing the random length of the  alternating environment (called $T$-stage and $S$-stage), and we assume a $T$-stage starts at time 0.    The model grows according to
 $\{\vk X^{(i)}(t), t\ge 0\}$ during the $i$th $T$-stage and according to $\{\vk Y^{(i)}(t), t\ge 0\}$ during the $i$th $S$-stage.

Based on the above we   define an alternating renewal process with renewal epochs
$$0=V_0< V_1<V_2<V_3<\cdots$$
with
$V_i=(T_1+S_1)+\cdots +(T_i+S_i)$ which is the $i$th environment cycle time.  Then the resulting $n$-dimensional process $\vk Z(t)=(Z_1(t),\cdots, Z_n(t))$, is defined as
\BQNY
\vk Z(t):=\left\{
  \begin{array}{ll}
 \vk Z(V_i)+\vk X^{(i+1)}(t-V_i), & \hbox{if \ $V_i <  t\le V_i+T_{i+1}$;} \\
\vk Z(V_i)+ \vk X^{(i+1)}(T_{i+1})+\vk Y^{(i+1)}(t-V_i-T_{i+1}), & \hbox{if \ $ V_i+T_{i+1} < t \le V_{i+1}$}.
  \end{array}
\right.
\EQNY
Note that this is a multi-dimensional regenerative process with regeneration epochs $V_i$.  This is a generalization of the one-dimensional  model discussed in \cite{KellaWhitt}. %It is worth mentioningn a very recent paper \cite{Jac19}, a different multi-dimensional regenerative process is introduced and studied.  ???

We assume that $\{\{\vk X^{(i)}(t), t\ge 0\}: i\ge 1\}$ and $\{\{\vk Y^{(i)}(t), t\ge 0\}: i\ge 1\}$ are independent samples of $\{\vk X(t), t\ge 0\}$ and $\{\vk Y(t), t\ge 0\}$, respectively, where
\BQNY
&&X_j(t)=B_{H_j}(t)+p_j t,\ \ \ t\ge0,\ \ \ \ 1\le j\le n,\\
&&Y_j(t)=\widetilde B_{\widetilde H_j} (t)-q_j t,\ \ \ t\ge0,\ \ \ \ 1\le j\le n,
\EQNY
with all the fBm's $B_{H_j}, \widetilde B_{\widetilde H_j}$  being mutually independent and $p_j, q_j>0, 1\le j\le n$. Suppose that $(T_i,S_i), i\ge 1$ are independent samples of $(T,S)$ and $T$ is regularly varying with index $\lambda>1$.
We further assume that  %for any $c>0$
\BQN\label{eq:ST}
\pk{S>x}=o\LT(\pk{T>x}\RT), \ \ \
%\EQN
 %and
%\BQN\label{eq:pq}
p_j \E{T} < q_j \E{S} <\IF\   \ 1\le j\le n.
\EQN
%for which it is necessary that $\lambda>1$.

For notational simplicity we shall restrict ourselves to the 2-dimensional case.
The general $n$-dimensional problem can be analysed similarly. Thus, for the rest of this section and related proofs in Section \ref{Sec:Qu}, all vectors (or multi-dimensional processes) are considered to be two-dimensional ones.

We are interested in the asymptotics of the following tail probability
\BQNY
Q(u):=\pk{\exists n \ge 1:\ \sup_{t\in [V_{n-1}, V_n]} Z_1(t) > a_1 u,  \sup_{s \in [V_{n-1}, V_n]} Z_2 (s)> a_2 u},\ \ \ u\to\IF,
\EQNY
with $a_1, a_2>0$. %Under the second assumption in \eqref{eq:ST}, we can ensure that $Q(u)$ is a rare event for large $u$.
In the fluid queueing context, $Q(u)$ can be interpreted as the probability that both buffers  overflow in some environment cycle.
In the insurance context, $Q(u)$ can be interpreted as the probability that in some business cycle the two lines of business of the insurer are both ruined (not necessarily at the same time). Similar one-dimensional models have been discussed   in the literature; see, e.g., \cite{ZBD05, PZ07, AsmAlb10}.

We introduce the following notation:
\BQN
&&\vk U^{(n)}=(U_1^{(n)}, U_2^{(n)}):=\vk Z(V_{n}) -\vk Z(V_{n-1}),\ \ n\ge 1,\ \ \  \vk U^{(0)}=\vk0, \label{eq:Un}\\
&&\vk M^{(n)}=(M_1^{(n)}, M_2^{(n)}):=\LT(\sup_{t\in [V_{n-1}, V_n)} Z_1(t)-Z_1(V_{n-1}), \sup_{s\in [V_{n-1}, V_n)} Z_2(s)-Z_2(V_{n-1})\RT), \ \ n\ge 1.\label{eq:Mn}
\EQN
Then we have
\BQNY
Q(u)=\pk{\exists n \ge 1:\ \sum_{i=1}^n  U^{(i-1)}_{1}+ M_1^{(n)} > a_1 u,  \sum_{i=1}^n  U^{(i-1)}_{2}+ M_2^{(n)} > a_2 u}.
\EQNY
Note that $\vk U^{(n)}, n\ge0$ and $\vk M^{(n)}, n\ge 0$ are both IID sequences. By the second assumption in \eqref{eq:ST} we have
\BQN\label{def-vkc}
 \E{\vk U^{(1)}}=(p_1\E{T}-q_1\E{S},p_2\E{T}-q_2\E{S})=:-\vk c <\vk 0,
\EQN
which ensures that %$Q(u)$
the event in the above probability is a rare event for large $u$, i.e., $Q(u)\to 0,$ as $u\to\IF$.

It is noted that our question now becomes an exit problem of a {\it 2-dimensional  perturbed random walk}. The exit problems of multi-dimensional random walk  has been discussed in many papers, e.g., \cite{Hult06}. %, where it is assumed that the increments  are regularly varying distributed .
However, it seems that multi-dimensional perturbed random walk has not been  discussed in the existing literature.

%Denote $ \vk c =-\E{\vk U^{(1)}}.$

%We have also to check (this follows from the tail equivalence Lemma \ref{lem:MU} and $\lambda>1$; maybe not used...)
%$$
%\E{\vk M^{(1)}}<\IF.
%$$

Since $T$ is regularly varying with index $\lambda>1$, we have that
\BQN\label{def:wtT}
\wtT:=(p_1 T, p_2 T)
\EQN
is regularly varying with index $\lambda$ and some limiting measure $\mu$ (whose form depends on the norm $| \cdot |$ that is chosen). % such that
%\BQNY
%\frac{\pk{x^{-1}\vk{\mathcal{T}} \in \cdot}}{\pk{\abs{\vk{\mathcal{T}}}>x}} \ \vto \ \mu(\cdot )\ \ \ \ \text{on}\ \mcB( [\vk 0,\vk \IF] \setminus \{\vk 0\}).
%\EQNY
%Note that the support of $\mu$ is on the line $\{\vk x=(x_1, x_2)\ge \vk0: x_1/p_1=x_2/p_2\}$. Particularly,
We present next the main result of this section, leaving its proof to Section \ref{Sec:Qu}.
\BT \label{Thm:Qu}
Under the above assumptions on regenerative model ${\vk Z}(t), t\ge0$, % and environment time lengths $T$ and $S$,
we have that, as $u\to\IF,$ %for any $\vk a>\vk 0$
\BQNY
Q(u)
\sim \int_0^\IF \mu( (v\vk c+\vk a, \vk \IF])  dv \  \pk{\abs{\wtT}>u} u,
\EQNY
where $\vk c$ and $\wtT$ is given by \eqref{def-vkc} and \eqref{def:wtT}, respectively.
\ET
\begin{remark}
Consider $\abs{\  \cdot\  }$ to be the $L^1$ norm in Theorem \ref{Thm:Qu}. %Therefore,  for any $\vk a\ge \vk 0, \vk a\neq \vk 0$,
We have
\BQNY\label{eq:mua}
\mu([\vk a, \vk \IF])=\LT( (p_1+p_2)\max(a_1/p_1, a_2/p_2)  \RT)^{-\lambda},
\EQNY
and thus, as $u\to\IF,$
\BQNY
Q(u) \sim
\int_0^\IF \max( (a_1+c_1 v)/p_1, (a_2+c_2 v)/p_2) ^{-\lambda} dv\   \pk{T> u} u.
\EQNY
\end{remark}

\section{Proof of main results} \label{Sec:promain}

This section is devoted to the proof of Theorem \ref{Thm:Main} followed by a short proof of Example \ref{Exm}.

First we give a result in line with  Proposition \ref{CorXc}.
Note that in the proof of the main results in \cite{dieker2005extremes} the minimum point $t_u^*$ of the function
\BQNY
f_u(t) :=\frac{u(1+t)}{\sigma(ut/c)}, \ \ t\ge 0,
\EQNY
plays an important role. It has been discussed therein that $t_u^*$ converges, as $u\to\IF,$ to $t^*:=H/(1-H)$ which is the unique minimum point of $\lim_{u\to\IF} f_u(t) \sigma(u)/u= (1+t)/ (t/c)^H, t\ge 0$. In this sense, $t_u^*$ is asymptotically unique.
We have the following corollary of \cite{dieker2005extremes}, which is  useful for the proofs below. %s, the proof of it is displayed in Section \ref{Sec:pro}.
\BEL \label{Prop:2}
Let $X(t), t\ge0$  be an a.s. continuous centered Gaussian process with stationary increments and $X(0)=0$. Suppose that {\bf C1--C4} hold. For any fixed $0<\vn<t^*/c$, we have, as $u\to\IF,$
\BQNY
\pk{\sup_{t\in 0, (t^*/c +\vn) u]} ( X(t)-ct) >u}\sim \psi(u), % ,\ \ \ \pk{\sup_{t\in[0, (t^*/c -\vn) u]} ( X(t)-ct) >u} =o(P_c(u)).
\EQNY
with $\psi(u)$ the same as  in  Proposition \ref{CorXc}.  Furthermore, we have that for any $\gamma>0$
\BQNY
\lim_{u\to\IF}\frac{\pk{  \sup_{t\in[0,(t^*/c-\vn) u ]} (X(t) -c t) >u    } }{\psi(u) u^{-\gamma}}=0.
\EQNY

\EEL

{\bf Proof of Lemma \ref{Prop:2}:} Note that
\BQNY
\pk{\sup_{t\in [0, (t^*/c +\vn) u]} ( X(t)-ct) >u}=
\pk{\sup_{t\in [0, (t^* + c\vn)]} \frac{X(ut/c)}{1+t}  >u}.
\EQNY
The first claim follows from \cite{dieker2005extremes}, as the main interval which determines the asymptotics is in$[0, (t^* + c\vn)]$ (see Lemma 7  and the comments in Section 2.1 therein). Similarly, we have
\BQNY
\pk{\sup_{t\in [0, (t^*/c -\vn) u]} ( X(t)-ct) >u}=
\pk{\sup_{t\in [0, (t^* - c\vn)]} \frac{X(ut/c)}{1+t}  >u}.
\EQNY
Since $t^*_u$ is asymptotically unique and   $\lim_{u\to\IF}t^*_u=t^*$, we can show that for all $u$ large
\BQNY
\inf_{t\in [0, (t^* - c\vn)]} f_u(t) \ge \rho f_u(t^*_u)=\rho \inf_{t\ge 0} f_u(t)
\EQNY
for some $\rho>1$. Thus, by similar arguments as in the proof of  Lemma 7 in \cite{dieker2005extremes} using the Borel inequality we conclude the second claim.
\QED

The following lemma is crucial for the proof of \netheo{Thm:Main}.

\BEL\label{Lem:XT}
%Let $\vk X(t), t\ge 0$ be a multi-dimensional centered Gaussian process with stationary increments that satisfies Assumption I, and
Let $X_i(t), t\ge 0$,  $i=1,2,\cdots, n_0 (< n)$ be independent centered Gaussian processes with stationary increments, and
let $\vk{\mathcal{T}}$ be an independent regularly varying random vector with index $\alpha$ and limiting measure $\nu$.   Suppose that all of
$\sigma_i(t), i=1,2,\cdots, n_0$ satisfy the assumptions {\bf C1}-{\bf C3} with the parameters involved indexed by $i$, which further satisfy that,
$\overleftarrow{\sigma_i}(u)\sim k_i \overleftarrow{\sigma}_1 (u) $ for some positive constants $k_i,i=1,2,\cdots,m\leq n_0$ and   $\overleftarrow{\sigma_{j} }(u)=o(\overleftarrow{\sigma_1 }(u))$ for all $j=m+1,\cdots, n_0$. Then, for any increasing to infinity functions $h_i(u), n_0+1\le i\le n$ such that $h_i(u)=o(\overleftarrow{\sigma_1}(u)), n_0+1\le i\le n$, and any $a_i>0$,
\BQNY
\pk{ \cap_{i=1}^{n_0} \LT( \sup_{t\in[0,\mathcal{T}_i]} X_{i}(t)>a_i u \RT), \cap_{i=n_0+1}^n \LT(\mathcal{T}_i>h_i(u) \RT)} \sim  \widetilde{\nu}((\vk k \vk a_{m,0}^{1/{\vk H}},\vk\IF]) \ \pk{\abs{\vk {\mathcal T}} > \isig1u},
%\pk{ \cap_{i=1}^{m} \LT( \mathcal{T}_i \xi_i^{1/H_i}  >b_i^{1/H_i} \cL{\overleftarrow{\sigma_i}}(u)  \RT) } .
\EQNY
where $\widetilde{\nu}$ is defined in \eqref{def:nutil} and $\vk k\vk a_{m,0}^{1/{\vk H}}=(k_1a_1^{1/H_1}, \cdots, k_ma_m^{1/H_m}, 0\cdots, 0)$ with $H_1=H_2=\cdots=H_m$.
\COM{$\widetilde{\nu}(K)=\E{\nu(\vk \xi^{-1/{\vk H}} K)}$, with $\xi_i$ defined in \eqref{eq:xieta} and $\vk \xi^{-1/{\vk H}} K=\{(\xi_1^{-1/{H_1}}d_1,\cdots, \xi_n^{-1/{H_n}}d_n), (d_1,\cdots,d_n)\in K\}$ for any $K\subset  [0,\IF]^n \setminus \{\vk 0\}$,}
\EEL

\prooflem{Lem:XT} We use a similar argument as in the proof of Theorem 2.1 in \cite{debicki2004supremum} to verify our conclusion.
For notational convenience denote
\BQNY
H(u)=:\pk{ \cap_{i=1}^{n_0} \LT( \sup_{t\in[0,\mathcal{T}_i]} X_{i}(t)>a_i u \RT), \cap_{i=n_0+1}^n \LT(\mathcal{T}_i>h_i(u) \RT)}.
\EQNY
We first give a asymptotically lower bound for $H(u)$. Let $G(\vk x)= \pk{\vk{\mathcal{T}} \le \vk x}$ be the distribution function of $\vk{\mathcal{T}}$. Note that, for any constants $0<r<R$,
\BQN\label{eq:HJJ}\nonumber
%&&\pk{ \cap_{i=1}^{n_0} \LT( \sup_{t\in[0,\mathcal{T}_i]} X_{i}(t)>b_i u \RT), \cap_{i=n_0+1}^n \LT(\mathcal{T}_i>h_i( u) \RT)} \\
H(u)&\geq&\pk{ \cap_{i=1}^{n_0} \LT( \sup_{t\in[0,\mathcal{T}_i]} X_{i}(t)>a_i u \RT), \cap_{i=1}^m (r\isig1u \le \mathcal{T}_i\le R\isig1u),\cap_{i=m+1}^n\LT(\MT_i>r\isig1u\RT)}\\ \nonumber
&=& \oint_{[r,R]^{m}\times (r,\infty)^{n-m}} \pk{\cap_{i=1}^{n_0}  \LT(\sup_{t\in[0, \overleftarrow{\sigma_1}(u)t_i]} X_{i}(t)>a_i u\RT)}
 dG(\overleftarrow{\sigma_1}(u) t_1,\cdots, \overleftarrow{\sigma_{1}}(u) t_{n})\\ \nonumber
 &=& \oint_{[r,R]^{m}\times (r,\infty)^{n-m}}\prod_{i=1}^{n_0}  \pk{ \sup_{s\in[0,1]} X_{i}^{u,t_i}(s)>a_i u_i(t_i)}
 dG(\overleftarrow{\sigma_1}(u) t_1,\cdots, \overleftarrow{\sigma_{1}}(u) t_{n})
\EQN
holds for sufficiently large $u$, where
\BQNY
X_{i}^{u,t_i}(s)=:\frac{X_i(\isig1u t_i s)}{\sigma_i(\isig1u t_i)},
u_i(t_i)=:\frac{u}{\sigma_i(\isig1u t_i)}, s\in[0,1],
(t_1,t_2,\cdots,t_{n_0})\in [r,R]^m\times(r,\infty)^{n_0-m}.
%t_i\in[b,B],i=1,\cdots,m,t_i\in(b,\infty),i=m+1,\cdots,n_0.
\EQNY
By Lemma 5.2 in \cite{debicki2004supremum}, we know that, as $u\to\infty$, the processes $X_{i}^{u,t_i}(s)$ converges weakly in $C([0,1])$ to $B_{H_i}(s)$, uniformly in $t_i\in(r,\infty)$, respectively for $i=1,2,\cdots,n_0$. Further, according to the assumptions on $\sigma_i(t)$, Theorem 1.5.2 and Theorem 1.5.6 in \cite{bingham1989regular}, we have, as $u\to\infty$, $u_i(t_i)$ converges to $k_i ^{H_i}t_i^{-H_i}$ uniformly in $t_i\in[r,R]$, respectively for $i=1,2,\cdots,m$, and $u_i(t_i)$ converges to $0$ uniformly in $t_i\in[r,\IF)$, respectively for $i=m+1,\cdots,n_0$. Then, by the continuous mapping theorem and recalling $\xi_i$ defined in \eqref{eq:xieta} is a continuous random variable (e.g., \cite{ZN03}), we get
\BQN\label{ieq:lbHu}
H(u)&\gtrsim& \oint_{[r,R]^{m}\times (r,\infty)^{n-m}}\prod_{i=1}^{m}  \pk{ \sup_{s\in[0,1]} B_{H_i}(s)>a_i k_i^{H_i}t_i^{-H_i}}
 dG(\overleftarrow{\sigma_1}(u) t_1,\cdots, \overleftarrow{\sigma_{1}}(u) t_{n})\\ \nonumber
&=& \pk{\cap_{i=1}^m\LT(\xi_i^{\frac{1}{H_i}}\MT_i>k_ia_i^{\frac{1}{H_i}}\isig1u\RT),
\cap_{i=1}^m\LT( r\isig1u\leq\MT_i\leq R\isig1u\RT),\cap_{i=m+1}^n\LT(\MT_i>r\isig1u\RT)}\\ \nonumber
%&=& \pk{A_1(u)} - \pk{A_1(u)\cap A_2(u)},
&=& J_1(u) -J_2(u),
%&=& \pk{\cap_{i=1}^m\LT(\xi^{\frac{1}{H_i}}\MT_i>(k_ia_i)^{\frac{1}{H_i}}\isig1u\RT),
%\cap_{i=m+1}^n\LT( r\isig1u\leq\MT_i\leq R\isig1u\RT)}\\
%&& -\ \pk{\cap_{i=1}^m\LT(\xi^{\frac{1}{H_i}}\MT_i>(k_ia_i)^{\1Hi}\isig1u\RT),
%\cap_{i=m+1}^n\LT( r\isig1u\leq\MT_i\leq R\isig1u\RT),
%\cup_{i=1}^{m} \LT( (\MT_i< r\isig1u)   \cup (\MT_i>R\isig1u )\RT)}
\EQN
where
\BQNY
%A_1(u)&=:&\cap_{i=1}^m\LT(\xi^{\frac{1}{H_i}}\MT_i>(k_ia_i)^{\frac{1}{H_i}}\isig1u\RT)
%\cap_{i=m+1}^n\LT(\MT_i>r\isig1u\RT),\\
%A_2(u)&=:&\cup_{i=1}^{m} \LT( (\MT_i< r\isig1u)   \cup (\MT_i>R\isig1u )\RT).
J_1(u)&=:&\pk{\cap_{i=1}^m\LT(\xi_i^{\frac{1}{H_i}}\MT_i>k_ia_i^{\frac{1}{H_i}}\isig1u\RT),
\cap_{i=m+1}^n\LT(\MT_i>r\isig1u\RT)},\\
J_2(u)&=:&\pk{\cap_{i=1}^m\LT(\xi_i^{\frac{1}{H_i}}\MT_i>k_ia_i^{\frac{1}{H_i}}\isig1u\RT),
\cap_{i=m+1}^n\LT(\MT_i>r\isig1u\RT),\cup_{i=1}^{m} \LT( (\MT_i< r\isig1u)   \cup (\MT_i>R\isig1u )\RT)}
\EQNY
Putting ${\vk \eta}=(\xi_1^{1/H_1},\cdots,\xi_m^{1/H_m},1,\cdots,1)$, then by \nelem{Lem:2} and the continuity of the limiting measure $\widehat\nu$ defined therein, we have
\BQN\label{lim:J1u}
\lim_{r\to 0}\lim_{u\to\infty}\frac{J_1(u)}{\pk{\abs{\vk \MT}>\isig1u}}
=\widetilde\nu((\vk k{\vk a}_{m,0}^{1/{\vk H}},\vk\infty]).
\EQN

Furthermore,
\BQNY
%&&\pk{ \cap_{i=1}^{n_0} \LT( \xi ^{1/H_i} \mathcal{T}_i>b_i^{1/H_i} \cL{\overleftarrow{\sigma_i}}(u)  \RT), \cap_{i=n_0+1}^n \LT(\mathcal{T}_i>b_i u \RT), \cup_{i=1}^{n_0} ( (\mathcal{T}_i< a\cL{\overleftarrow{\sigma_i}}(u))   \cup (\mathcal{T}_i>A\cL{\overleftarrow{\sigma_i}}(u) )) }\\
J_2(u)\leq  \sum_{i=1}^{m} \LT( \pk{ \xi_i ^{\1Hi} \MT_i>k_i a_i^{\1Hi} \isig1u, \mathcal{T}_i< r\isig1u} +\pk{\mathcal{T}_i>R\isig1u}\RT).
\EQNY
Then, by the fact that $\abs{\vk{\mathcal{T}}}$  is regularly varying with index $\alpha$, and using the same arguments as in the the proof of Theorem 2.1 in \cite{debicki2004supremum} (see the asymptotic for integral $I_4$ and (5.14) therein), we conclude that
\BQN\label{lim:J2u}
\lim_{r\to0,R\to\infty} \limsup_{u\to\IF}\frac{J_2(u)}{ \pk{\abs{\vk{\mathcal{T}}} > \overleftarrow{\sigma_1 }(u)}} =0,
\EQN
which combined with \eqref{ieq:lbHu} and \eqref{lim:J1u} yields
\BQN\label{lb:Hu}
\lim_{r\to 0,R\to\IF}\liminf_{u\to\infty}\frac{H(u)}{\pk{\abs{\vk \MT}>\isig1u}}
\geq\widetilde\nu((\vk k{\vk a}_{m,0}^{1/{\vk H}},\vk\infty]).
\EQN
Next, we give an asymptotic upper bound for $H(u)$. Note
\BQNY
H(u)&\le &%\pk{ \cap_{i=1}^{n_0} \LT( \sup_{t\in[0,\mathcal{T}_i]} X_{i}(t)>b_i u \RT), \cap_{i=n_0+1}^n \LT(\mathcal{T}_i>h_i( u) \RT)} \le
\pk{ \cap_{i=1}^{m} \LT( \sup_{t\in[0,\mathcal{T}_i]} X_{i}(t)>a_i u \RT) } \\
&=&   \pk{ \cap_{i=1}^{m} \LT( \sup_{t\in[0,\mathcal{T}_i]} X_{i}(t)>a_i u \RT), \cap_{i=1}^m (r\isig1u \le \mathcal{T}_i\le R\isig1u)}\\
&&+\ \ \pk{ \cap_{i=1}^{m} \LT( \sup_{t\in[0,\mathcal{T}_i]} X_{i}(t)>a_i u \RT), \cup_{i=1}^{m} \LT( (\MT_i< r\isig1u)   \cup (\MT_i>R\isig1u )\RT)}\\
&=:& J_3(u)+J_4(u).
\EQNY
By the same reasoning as that used in the deduction for \eqref{ieq:lbHu}, we can show that
\COM{that used in the deduction of \eqref{ieq:lbHu}, we have
\BQNY
J_3(u)&\sim& \pk{ \cap_{i=1}^m\LT(\xi^{\frac{1}{H_i}}\MT_i>k_ia_i^{\frac{1}{H_i}}\isig1u\RT),
\cap_{i=1}^m (r\isig1u \le \mathcal{T}_i\le R\isig1u)}\\
&=& \pk{ \cap_{i=1}^m\LT(\xi^{\frac{1}{H_i}}\MT_i>(k_ia_i)^{\frac{1}{H_i}}\isig1u\RT)}\\
&&-\ \ \pk{\cap_{i=1}^m\LT(\xi^{\frac{1}{H_i}}\MT_i>(k_ia_i)^{\frac{1}{H_i}}\isig1u\RT)
,\cup_{i=1}^{m} \LT( (\MT_i< r\isig1u)   \cup (\MT_i>R\isig1u )\RT)}\\
&=:&J_5(u)-J_6(u).
\EQNY
Then, applying similar analysis as used in \eqref{lim:J1u} and \eqref{lim:J2u} yields
\BQNY
\lim_{u\to\infty}\frac{J_5(u)}{\pk{\abs{\vk \MT}>\isig1u}}
=\widetilde\nu(((\vk k{\vk a})^{1/{\vk H}},\vk\infty])\ \ \ \textrm{and} \
\lim_{r\to0,R\to\infty} \limsup_{u\to\IF}\frac{J_6(u)}{ \pk{\abs{\vk{\mathcal{T}}} > \overleftarrow{\sigma_1 }(u)}} =0,
\EQNY
and thus}
\BQN\label{lim:J3u}
\lim_{r\to0,R\to\infty} \lim_{u\to\infty}\frac{J_3(u)}{\pk{\abs{\vk \MT}>\isig1u}}
=\widetilde\nu((\vk k{\vk a}_{m,0}^{1/{\vk H}},\vk\infty]).
\EQN
Moreover,
\BQNY
J_4(u)\leq \sum_{i=1}^{m} \LT( \pk{ \sup_{t\in[0,\mathcal{T}_i]} X_{i}(t)>a_i u, \mathcal{T}_i< r\isig1u} +\pk{\mathcal{T}_i>R\isig1u} \RT).
\EQNY
Thus, by the same arguments as in the proof of Theorem 2.1 in \cite{debicki2004supremum} (see the asymptotics for integrals $I_1, I_2, I_4$ therein), we conclude
\BQNY
\lim_{r\to0,R\to\infty} \limsup_{u\to\IF}\frac{J_4(u)}{ \pk{\abs{\vk{\mathcal{T}}} > \overleftarrow{\sigma_1 }(u)}} =0,
\EQNY
which together with \eqref{lim:J3u} implies that
\BQN\label{ub:Hu}
\lim_{r\to 0,R\to\IF}\limsup_{u\to\infty}\frac{H(u)}{\pk{\abs{\vk \MT}>\isig1u}}
\leq\widetilde\nu((\vk k{\vk a}_{m,0}^{1/{\vk H}},\vk\infty]).
\EQN
Notice that by the assumptions on $\{\overleftarrow{\sigma_i}(u)\}_{i=1}^{m}$, we in fact have $H_1=H_2=\cdots=H_m$. Consequently, combing \eqref{lb:Hu} and \eqref{ub:Hu} we  complete  the proof. \QED

\prooftheo{Thm:Main} %Before the beginning of proof,
We use in the following the convention that $\cap_{i=1}^0=\Omega$, the sample space. We first verify the claim for case (i),  $n_0>0$. For arbitrarily small $\vn>0$, we have
\BQNY
%&&\pk{\cap_{i=1}^n \LT(\sup_{t\in[0,\mathcal{T}_i]} (X_{i}(t) +c_i t) >a_i u \RT)}\\
P(u)&\ge&  \mathbb{P}\Bigg\{\cap_{i=1}^{n_-} \LT(\sup_{t\in[0,\mathcal{T}_i]} (X_{i}(t) +c_i t) >a_i u,\mathcal{T}_i>(t^*_i/\abs{c_i}+\vn) u \RT), \cap_{i=n_-+1}^{n_-+n_0} \LT( \sup_{t\in[0,\mathcal{T}_i]} X_{i}(t) >a_i u \RT),  \\
&&\ \ \ \ \  \cap_{i=n_-+n_0+1}^{n} \LT(\sup_{t\in[0,\mathcal{T}_i]} (X_{i}(t) +c_i t) >a_i u,\MT_i>\frac{a_i+\vn}{c_i}u \RT)\Bigg\}\\
&\geq&  \mathbb{P}\Bigg\{\cap_{i=1}^{n_-} \LT(\sup_{t\in[0, (t^*_i/\abs{c_i}+\vn) u ]} (X_{i}(t) +c_i t) >a_i u,\mathcal{T}_i>(t^*_i/\abs{c_i}+\vn) u \RT),  \\
&&\ \ \ \ \ \cap_{i=n_-+1}^{n_-+n_0} \LT( \sup_{t\in[0,\mathcal{T}_i]} X_{i}(t) >a_i u \RT), \cap_{i=n_-+n_0+1}^{n} \LT( X_{i}\LT(\frac{a_i+\vn}{c_i}u\RT) >-\vn u,\MT_i>\frac{a_i+\vn}{c_i}u \RT)\Bigg\}\\
&=& Q_1(u)\times Q_2(u)\times Q_3(u),
\EQNY
where
\BQNY
Q_1(u)&:=&\pk{\cap_{i=1}^{n_-} \LT(\sup_{t\in[0, (t^*_i/\abs{c_i}+\vn) u ]} X_{i}(t) +c_i t >a_i u \RT)}\\
Q_2(u)&:=&\pk{\cap_{i=1}^{n_-} \LT( \mathcal{T}_i>(t^*_i/\abs{c_i}+\vn) u \RT), \cap_{i=n_-+1}^{n_-+n_0} \LT( \sup_{t\in[0,\mathcal{T}_i]} X_{i}(t)>a_i u \RT), \cap_{i=n_-+n_0+1}^n \LT( \MT_i>\frac{a_i+\vn}{c_i}u \RT)},\\
Q_3(u)&:=&\prod_{i=n_-+n_0+1}^n \pk{N_i> \frac{-\vn u}{ \sigma_i(\frac{a_i+\vn}{c_i}u)} }\to1,\ \ \ u\to\IF,
\EQNY
with $N_i,i=n_-+n_0+1,\cdots,n$ being standard Normal distributed random variables. %, independent of $\vk X$ and $\vk \MT$.
By \nelem{Prop:2}, we know, as $u\to\IF$,
\BQNY
Q_1(u)\sim \prod_{i=1}^{n_-}\psi_i(a_i u).
\EQNY
Further, according to the assumptions on $\sigma_i$'s and \nelem{Lem:XT}, we get
\BQNY
\lim_{\vn\to0}\lim_{u\to\IF}\frac{Q_2(u)}{\pk{\abs{\vk \MT}>\overleftarrow{\sigma}_{n_-+1}(u)}}=\widetilde{\nu}((\vk k\vk a_0^{1/{H_{n_-+1}}},\vk\IF]), %\quad \textrm{and}\quad \lim_{u\to\IF}Q_3(u)=1,
\EQNY
and thus
\BQNY
P(u)\gtrsim \widetilde{\nu}((\vk k\vk a_0^{1/{H_{n_-+1}}},\vk\IF]) \pk{\abs{\vk \MT}>\overleftarrow{\sigma}_{n_-+1}(u)} \prod_{i=1}^{n_-}\psi_i(a_i u),\quad u\to\IF.
\EQNY
%Moreover, using \neprop{Prop:2} and \nelem{Lem:XT} again, we have
Similarly, we can show
\BQNY
P(u)&\le & \pk{\cap_{i=1}^{n_-} \LT(\sup_{t\in[0,\IF)} X_{i}(t) +c_i t >a_i u \RT), \cap_{i=n_-+1}^{n_-+n_0}  \LT(\sup_{t\in[0,\mathcal{T}_i]} X_{i}(t)>a_i u \RT)} \\
&&\sim\ \widetilde{\nu}((\vk k\vk a_0^{1/{H_{n_-+1}}},\vk\IF]) \pk{\abs{\vk \MT}>\overleftarrow{\sigma}_{n_-+1}(u)} \prod_{i=1}^{n_-}\psi_i(a_i u),\quad u\to\IF.
\EQNY
This completes the proof of case (i).

Next we consider case (ii), $n_0=0$.  Similarly as in case (i) we have, for any small $\vn>0$
 \BQNY
P(u) %&\ge&  \mathbb{P}\Bigg\{\cap_{i=1}^{n_-} \LT(\sup_{t\in[0,\mathcal{T}_i]} (X_{i}(t) +c_i t) >a_i u,\mathcal{T}_i>(t^*_i/\abs{c_i}+\vn) u \RT), \\
%&&\ \ \ \ \  \cap_{i=n_-+1}^{n} \LT(\sup_{t\in[0,\mathcal{T}_i]} (X_{i}(t) +c_i t) >a_i u,\MT_i>\frac{a_i+\vn}{c_i}u \RT)\Bigg\}\\
&\geq&  \mathbb{P}\Bigg\{\cap_{i=1}^{n_-} \LT(\sup_{t\in[0, (t^*_i/\abs{c_i}+\vn) u ]} (X_{i}(t) +c_i t) >a_i u,\mathcal{T}_i>(t^*_i/\abs{c_i}+\vn) u \RT),  \\
&&\ \ \ \ \ \cap_{i=n_-+1}^{n} \LT( X_{i}\LT(\frac{a_i+\vn}{c_i}u\RT) >-\vn u,\MT_i>\frac{a_i+\vn}{c_i}u \RT)\Bigg\}\\
&=& Q_1(u)\times Q_3(u)\times Q_4(u),
\EQNY
where
\BQNY
Q_4(u):=\pk{\cap_{i=1}^{n_-} \LT(\mathcal{T}_i>(t^*_i/\abs{c_i}+\vn) u \RT),
\cap_{i=n_-+1}^{n} \LT(\MT_i>\frac{a_i+\vn}{c_i}u \RT) }.
\EQNY
By \nelem{Lem:1}, we know
\BQNY
\lim_{\vn\to0}\lim_{u\to\IF}\frac{Q_4(u)}{\pk{\abs{\vk \MT}>u}}= \nu(\vk a_1, \vk\IF],
\EQNY
and thus
\BQNY
P(u)\gtrsim \nu(\vk a_1, \vk\IF] \pk{\abs{\vk \MT}>u} \prod_{i=1}^{n_-}\psi_i(a_i u),\quad u\to\IF.
\EQNY
For the upper bound, we have for any small  $\vn>0$
\BQNY
P(u)\le I_1(u)+I_2(u),
\EQNY
with
\BQNY
%&&\pk{\cap_{i=1}^n \LT(\sup_{t\in[0,\mathcal{T}_i]} X_{i}(t) +c_i t >a_i u \RT)}\\
I_1(u)&:=&\mathbb{P}\Bigg\{\cap_{i=1}^{n_-} \LT(\sup_{t\in[0,\mathcal{T}_i]} X_{i}(t) +c_i t >a_i u \RT), \cap_{i=1}^{n_-} \LT( \mathcal{T}_i>(t^*_i/\abs{c_i}-\vn) u \RT),   \cap_{i=n_-+1}^n \LT( \sup_{t\in[0,\mathcal{T}_i]} X_{i}(t) +c_i \mathcal{T}_i>a_i u \RT)\Bigg\},\\
I_2(u)&:=&    \mathbb{P}\Bigg\{\cap_{i=1}^{n_-} \LT(\sup_{t\in[0,\mathcal{T}_i]} X_{i}(t) +c_i t >a_i u \RT), \cup_{i=1}^{n_-} \LT( \mathcal{T}_i\le(t^*_i/\abs{c_i}-\vn) u \RT),   \cap_{i=n_-+1}^n \LT( \sup_{t\in[0,\mathcal{T}_i]} X_{i}(t)+c_i \mathcal{T}_i>a_i u \RT)\Bigg\}.%\\
%&&=: I_1(u)+I_2(u).
\EQNY
It follows that
\BQNY
I_1(u)&\le& \mathbb{P}\Bigg\{\cap_{i=1}^{n_-} \LT(\sup_{t\in[0,\IF)} X_{i}(t) +c_i t >a_i u \RT), \cap_{i=1}^{n_-} \LT( \mathcal{T}_i>(t^*_i/\abs{c_i}-\vn) u \RT),   \cap_{i=n_-+1}^n \LT( \sup_{t\in[0,\mathcal{T}_i]} X_{i}(t)+c_i \mathcal{T}_i>a_i u \RT)\Bigg\}\\
&=& \prod_{i=1}^{n_-} \psi_i(a_iu)    \pk{\cap_{i=1}^{n_-} \LT( \mathcal{T}_i>(t^*_i/\abs{c_i}-\vn) u \RT),   \cap_{i=n_-+1}^n \LT( \sup_{t\in[0,\mathcal{T}_i]} X_{i}(t)+c_i \mathcal{T}_i>a_i u \RT)}.
%\pk{\cap_{i=1}^{n_-} \LT( \mathcal{T}_i> t^*_i  u \RT), \cap_{i=n_-+1}^n \LT( c_i \mathcal{T}_i>a_i u \RT)}.
\EQNY
Next, we have for the small chosen $\vn>0$ %consider $i=n_-+1,\cdots,n$.
\BQNY
&& \pk{\cap_{i=1}^{n_-} \LT( \mathcal{T}_i>(t^*_i/\abs{c_i}-\vn) u \RT),   \cap_{i=n_-+1}^n \LT( \sup_{t\in[0,\mathcal{T}_i]} X_{i}(t) +c_i \mathcal{T}_i>a_i u \RT)}\\
&&= \pk{\cap_{i=1}^{n_-} \LT( \mathcal{T}_i>(t^*_i/\abs{c_i}-\vn) u \RT),   \cap_{i=n_-+1}^n \LT( \sup_{t\in[0,\mathcal{T}_i]} X_{i}(t) +c_i \mathcal{T}_i>a_i u,  \sup_{t\in[0,\mathcal{T}_i]} X_{i}(t) \le \vn u\RT)}\\
&&\ \ +\pk{\cap_{i=1}^{n_-} \LT( \mathcal{T}_i>(t^*_i/\abs{c_i}-\vn) u \RT),   \cap_{i=n_-+1}^n \LT( \sup_{t\in[0,\mathcal{T}_i]} X_{i}(t)+c_i \mathcal{T}_i>a_i u\RT), \cup_{i=n_-+1}^n \LT(\sup_{t\in[0,\mathcal{T}_i]} X_{i}(t) > \vn u\RT)}\\
&&\le \pk{\cap_{i=1}^{n_-} \LT( \mathcal{T}_i>(t^*_i/\abs{c_i}-\vn) u \RT),   \cap_{i=n_-+1}^n \LT(  c_i \mathcal{T}_i>(a_i-\vn) u \RT)} +\sum_{i=n_-+1}^n  \pk{  \sup_{t\in[0,\mathcal{T}_i]} X_{i}(t)> \vn u }. %\\
%&&=\pk{\cap_{i=1}^{n_-} \LT( \mathcal{T}_i>(t^*_i/\abs{c_i}-\vn) u \RT),   \cap_{i=n_-+1}^n \LT(  c_i \mathcal{T}_i>(a_i-\vn) u \RT)} (1+o(1)).
\EQNY
Furthermore, it follows from Theorem 2.1 in \cite{debicki2004supremum} that for any $i=n_-+1,\cdots, n$
\BQNY
\pk{  \sup_{t\in[0,\mathcal{T}_i]} X_{i}(t)> \vn u } \sim C_i(\vn) \pk{\mathcal{T}_i>\cL{\overleftarrow{\sigma_i}}(u)},\ \ \ \ u\to\IF,
\EQNY
with some constant $C_i(\vn)>0$. This implies that
\BQNY
\sum_{i=n_-+1}^n  \pk{  \sup_{t\in[0,\mathcal{T}_i]} X_{i}(t)> \vn u } =o(\pk{\abs{\vk{\mathcal{T}}}>u}),\quad u\to\IF.
\EQNY
Consequently, applying \nelem{Lem:1} and   letting $\vn \to 0$ we can obtain the required asymptotic upper bound, %given that
 if we can further show
\BQN\label{eq:I2u}
\lim_{u\to\IF}\frac{I_2(u)}{\prod_{i=1}^{n_-} \psi_i(a_iu)  \pk{\abs{\vk{\mathcal{T}}}>u }}=0.
\EQN
Indeed, we have
\BQN\label{eq:I_2}
I_2(u)&\le& \sum_{i=1}^{n_-} \mathbb{P}\Bigg\{\cap_{j=1}^{n_-} \LT(\sup_{t\in[0,\mathcal{T}_j]} X_{j}(t) +c_j t >a_j u \RT),    \mathcal{T}_i\le(t^*_i/\abs{c_i}-\vn) u   \Bigg\}\nonumber \\
&\le& \sum_{i=1}^{n_-}  \underset{j\neq i}{\prod_{j=1}^{n_-}} \psi_j(a_ju) \pk{  \sup_{t\in[0,(t^*_i/\abs{c_i}-\vn) u ]} X_{i}(t) +c_i t >a_i u    }.
\EQN
Furthermore, by \nelem{Prop:2} we have that for any $\gamma>0$
\BQNY
\lim_{u\to\IF}\frac{\pk{  \sup_{t\in[0,(t^*_i/\abs{c_i}-\vn) u ]} X_{i}(t) +c_i t >a_i u    } }{\psi_i(a_iu) u^{-\gamma}}=0,\quad i=1,2,\cdots,n_-,
\EQNY
which together with \eqref{eq:I_2} implies \eqref{eq:I2u}. This completes the proof. \QED

 {\bf Proof of Example \ref{Exm}: } The proof is based on the following obvious bounds
  \BQN\label{eq:ULB}
P_L(u)&:=&\pk{\cap_{i=1}^n \LT( ( B_{i}(Y_i(T)) +c_i Y_i(T) )> a_iu \RT)}\le P_B(u)\nonumber \\
&&  \le \pk{\cap_{i=1}^n \LT(\sup_{t\in[0, Y_i(T)]} ( B_{i}(t) +c_i t )> a_i u \RT)}=:P_U(u).
 \EQN
Since $\alpha_0<\min_{i=1}^n \alpha_i,$ by Lemma \ref{Lem:sum} we have that $\vk Y(T)$ is a multivariate regularly varying random vector with index $\alpha_0$ and the same limiting measure $\nu$ as that of $\vk S_0(T):=(S_0(T), \cdots, S_0(T))\in \R^n$, and further $\pk{\abs{\vk Y(T)}>x}\sim \pk{\abs{\vk S_0(T)}>x}, x\to\IF.$
% (different for different choice of norm). %satisfying $\nu([\vk a, \vk\IF])= (n \max(a_1,\cdots, a_n) )^{-\alpha}$ for any $\vk a \ge \vk0$ and $\vk a \neq \vk 0$.
The asymptotics of $P_U(u)$ can be obtained by applying Theorem \ref{Thm:Main}. Below we focus on $P_L(u)$.

 First, consider case (i) where $c_i>0$ for all $i=1,\cdots,n$. We have
 \BQNY
P_L(u) = \pk{\cap_{i=1}^n \LT( ( B_{i}(1 ) \sqrt{Y_i(T)}+c_i Y_i(T) )> a_iu \RT)}.
 \EQNY
 Thus, by \nelem{Lem:sum}   we obtain
 \BQNY
 P_L(u) \sim \pk{\cap_{i=1}^n \LT(  c_i S_0(T) > a_iu \RT)}\sim C_{\alpha_0,T} (\max_{i=1}^n (a_i/c_i)u)^{-\alpha_0}, \ \ u\to\IF,
 \EQNY
which is the same as the asymptotic upper bound obtained by using (ii) of Theorem \ref{Thm:Main}.

Next,   consider case (ii) where $c_i=0$ for all $i=1,\cdots,n$. We have
 \BQNY
P_L(u) = \pk{\cap_{i=1}^n \LT(  B_{i}(1 )  \sqrt{Y_i(T)} > a_iu \RT)} = \frac{1}{2^n} \pk{\cap_{i=1}^n \LT(  B_{i}(1 ) ^2  Y_i(T)  > (a_iu)^2 \RT)} .
 \EQNY
 Thus, by \nelem{Lem:2} and \nelem{Lem:sum}   we obtain
 \BQNY
 P_L(u)\asymp O(u^{-2\alpha_0}), \ \ u\to\IF,
 \EQNY
 which is the same as the asymptotic upper bound obtained by using (i) of  Theorem \ref{Thm:Main}.

 Finally, consider the case (iii) where $c_i<0$ for all $i=1,\cdots,n$. We have
   \BQNY
%&&\pk{\cap_{i=1}^n \LT( ( B_{i}(Y_i(T)) +c_i  Y_i(T) )> a_iu \RT)} \\
&&P_L(u)\ge \pk{\cap_{i=1}^n \LT(   B_{i}(Y_i(T)) +c_i  Y_i(T)>  a_iu, Y_i(T)\in [a_iu/\abs{c_i} -\sqrt u,  a_iu/\abs{c_i} +\sqrt u]\RT) } \\
&&\ge  \prod_{i=1}^n \LT(\min_{t\in [a_iu/\abs{c_i}-\sqrt u,  a_iu/\abs{c_i}+\sqrt u]} \pk{  B_{1}(t) +c_i t  > a_i u } \RT)  \pk{\cap_{i= 1}^n \LT( Y_i(T)\in [a_iu/\abs{c_i} -\sqrt u,  a_iu/\abs{c_i} +\sqrt u] \RT) }. %\\
%&\ge& \min_{t\in [1-u^{-1/2}, 1+u^{-1/2}]} \pk{  B_{1}(1)  >\frac{1+t}{\sqrt t}  \sqrt u }^{n_-}\\
%&&\times \pk{B_1(1)>-\sqrt {\frac{u}{u-\sqrt u}} }^{n_+} \pk{ S(T)\in [u-\sqrt u,  u+\sqrt u] }.
 \EQNY
 %Recall for the standard normal random variable $N_1$ we have %(see, e.g., \cite{AdlerTaylor})
 %\BQNY
% \pk{N_1>x}\sim \frac{1}{\sqrt{2\pi} x}  e^{-\frac{x^2}{2}}, \ \ x\to\IF.
 %\EQNY
%Letting $x_i(t,u)=(a_i-c_it)\sqrt u/\sqrt{t}$ and using the self-similarity of Brownian motion,
Recalling \eqref{eq:Psiu}, we derive that
\BQNY
\min_{t\in [a_iu/\abs{c_i}-\sqrt u,  a_iu/\abs{c_i}+\sqrt u]} \pk{  B_{1}(t) +c_i t  > a_i u }&=& \min_{t\in [a_i /\abs{c_i}-1/\sqrt u,  a_i /\abs{c_i}+1/\sqrt u]} \pk{  B_{1}(1)   > (a_i-c_it)\sqrt u/\sqrt{t}}\\
& \gtrsim & constant \cdot \frac{1}{\sqrt u} e^{2a_ic_i u+o(u)},\ \ \ u\to\IF.
\EQNY
Furthermore,
\BQN\nonumber
&&\pk{\cap_{i= 1}^n \LT( Y_i(T)\in [a_iu/\abs{c_i} -\sqrt u,  a_iu/\abs{c_i} +\sqrt u] \RT) }\\ \label{proSi}
&&\ge\prod_{i= 0}^n \pk{  S_i(T)\in [a_iu/\abs{2c_i} -\sqrt u/2,  a_iu/\abs{2c_i} +\sqrt u/2]  }.
\EQN
Due to the assumptions on the density functions of $S_i(T),i=0,1,\cdots,n$, then by Monotone Density Theorem (see e.g. in \cite{Mik99}),
we know that \eqref{proSi} is asymptotically larger than $C u^{-\beta}$ for some constants $C, \beta>0$.
%and the right hand side is a regularly varying function of $u$, due to the assumption that the density function of $S_i(T)$ is ultimately monotone  for all $i=0, 1,\cdots, n.$
Therefore,
$$
\ln P_L(u)   \gtrsim  2 \sum_{i=1}^n (a_ic_i) u,\ \ \ \ u\to\IF.
$$
The same asymptotic upper bound can be obtained by % Theorem \ref{Thm:Main} and
the fact that $\pk{\sup_{t>0}(B_i(t)+c_it)>a_iu}=e^{2a_ic_iu}$ for $c_i<0$. This completes the proof.
\QED

 %from Theorem \ref{Thm:Main} and the following obvious bounds

 \section{Proof of \netheo{Thm:Qu}} \label{Sec:Qu}

We first show one lemma which is  crucial for the proof of \netheo{Thm:Qu}.
%The proofs of them are displayed in the Appendix.

\BEL \label{lem:MU}
Let $\vk U^{(1)}$, $\vk M^{(1)}$ and $\wtT$ be given by \eqref{eq:Un}, \eqref{eq:Mn} and \eqref{def:wtT} respectively. Then, $\vk U^{(1)}, \vk M^{(1)}$ are both regularly varying with the same index $\lambda$ and limiting measure $\mu$ as that of $\wtT$. Moreover,
%Both  $\vk M^{(1)}$ and $\vk U^{(1)}$ are multivariate regularly varying random variables.
\BQNY
\pk{\abs{\vk U^{(1)}}> x} \sim \pk{\abs{\vk M^{(1)}} >x}\sim  \pk{\abs{\wtT}>x}, \ x\to \IF.
\EQNY
 \EEL

\prooflem{lem:MU} % We first consider $\vk U^{(1)}$.
First note that by self-similarity of fBm's
\BQNY
\vk U^{(1)}%\overset{D}=(B_{H_1}(T)+p_1 T+ \wtB_{\wtH_1}(S)-q_1 S,\ \  B_{H_2}(T)+p_2 T+ \wtB_{\wtH_2}(S)-q_2 S)
=(X_1^{(1)}(T_1)+Y_1^{(1)}(S_1), X_2^{(1)}(T_1)+Y_2^{(1)}(S_1))\overset{D}=(\vk \wtT+\vk Z_1+\vk Z_2+\vk Z_3),
\EQNY
where
%\BQNY
%&&{\vk Z}_1 =(p_1 T, p_2 T),\ \vk Z_2=(B_{H_1}(1) T^{H_1}, B_{H_2}(1) T^{H_2}),\\
%&&\vk Z_3=(\wtB_{\wtH_1}(1) S^{\wtH_1}, \wtB_{\wtH_2}(1) S^{\wtH_2}),\ \vk Z_4 =(-q_1 S, -q_2 S).
%\EQNY
\BQNY
&&\vk Z_1=(B_{H_1}(1) T^{H_1}, B_{H_2}(1) T^{H_2}),\
\vk Z_2=(\wtB_{\wtH_1}(1) S^{\wtH_1}, \wtB_{\wtH_2}(1) S^{\wtH_2}),\ \vk Z_3 =(-q_1 S, -q_2 S).
\EQNY
%Clearly, $\vk Z_1=\wtT$ is regularly varying with index $\lambda$ and limiting measure $\mu$ defined in \eqref{eq:mua}.
Since every two norms on $R^d$ are equivalent, then by the fact that $H_i, \wtH_i<1$ for $i=1,2$ and \eqref{eq:ST}, we have
  \BQNY
\max\Big(\pk{\abs{(T^{H_1},T^{H_2})}>x},\pk{\abs{(S^{\wtH_1},S^{\wtH_2})}>x},
\pk{\abs{\vk Z_3}>x}\Big)=o\Big(\pk{\abs{\wtT}>x}\Big),\ x\to\infty.
\EQNY
Thus, the claim for $\vk U^{(1)}$ follows directly by \nelem{Lem:sum}.
% if we can show that, as $x\to\IF,$
%\BQN
%&&\pk{\abs{\vk Z_2}>x} = o(\pk{\abs{\vk{\mathcal{T}}}>x} ), \label{eq:Z2}\\
%&& \pk{\abs{\vk Z_3}>x} = o(\pk{\abs{\vk{\mathcal{T}}}>x} ), \label{eq:Z3}\\
%&& \pk{\abs{\vk Z_4}>x} = o(\pk{\abs{\vk{\mathcal{T}}}>x} ).\label{eq:Z4}
%\EQN
%Note that \eqref{eq:Z4}  follows directly from \eqref{eq:ST}. Next we have  % from the well-known Breiman's theorem, as $x\to\IF$
%\BQNY
%\pk{\abs{\vk Z_2}>x} &\le& \pk{  \abs{B_{H_1}(1)} T^{H_1}+ \abs{  B_{H_2}(1)} T^{H_2}  >x}\\
%&\le&\pk{  2\abs{B_{H_1}(1)} T^{H_1}   >x} + \pk{  2\abs{B_{H_2}(1)} T^{H_2}   >x}  \\
%&=& o(\pk{  T >x} ),
%\EQNY
%where the last asymptotic equality follows from Lemma \ref{Lem:xiST} and the fact that $\pk{T^{H_i}   >x}=o(\pk{T >x}), i=1,2$.
% Similarly,  by \eqref{eq:ST} we have %as $x\to\IF,$
%\BQNY
%\pk{\abs{\vk Z_3}>x} %&\le&\pk{  2\abs{\wtB_{\wtH_1}(1)}   S^{ \wtH_1 }>x}+ \pk{  2\abs{\wtB_{\wtH_2}(1)}   S^{ \wtH_2 }>x}\\
%& \le& \pk{  \LT(2\abs{\wtB_{\wtH_1}(1)}\RT) ^{ 1/\wtH_1 }  S>x}+ \pk{  \LT(2\abs{\wtB_{\wtH_2}(1)} \RT)^{ 1/\wtH_2 }  S>x} \\
%& =&o(\pk{  T >x} ).
%\EQNY
%Consequently, since $\pk{\abs{\vk{\mathcal{T}}}>x}$ is regularly varying with index $\lambda,$ \eqref{eq:Z2} and \eqref{eq:Z3} are established. Thus, the  claim for $\vk U^{(1)}$ is established.

%Below, we focus on $\vk M^{(1)}$. We shall first show that, as $x\to\IF,$
%\BQN \label{eqLM1}
%\pk{ \abs{ \vk M^{(1)}} >x} \ \sim \ \pk{ \abs{\vk{\mathcal{T}}} >x}.
%\EQN
Next, note that
\BQNY
\vk M^{(1)}&\overset{D}=&\Big(\sup_{0\le t\le T+S} \LT(X_1(t) I_{(0\le t<T)} + (X_1(T) +Y_1(t-T)   ) I_{(T\le t<T+S)}\RT),\\
&&\ \ \ \ \ \ \ \ \sup_{0\le t\le T+S} \LT(X_2(t) I_{(0\le t<T)} + (X_2(T) +Y_2(t-T)   ) I_{(T\le t<T+S)} \RT)    \Big)=:\vk M,
\EQNY
then
\BQNY
\vk M\geq (X_1(T),X_2(T))\overset{D}=\wtT+\vk Z_1
\EQNY
and
\BQNY
\vk M&\leq& \Big(\sup_{0\leq t\leq T}B_{H_1}(t)+p_1T+\sup_{t\geq0}Y_1(t),
\sup_{0\leq t\leq T}B_{H_2}(t)+p_1T+\sup_{t\geq0}Y_2(t)\Big)\\
&\overset{D}=&
(\xi_1 T^{H_1} + \sup_{t\ge 0} Y_1(t),\xi_2 T^{H_2} +\sup_{t\ge 0} Y_2(t))+\wtT,
\EQNY
with $\xi_i$ defined in \eqref{eq:xieta}. By Corollary \ref{eq:psi_i}, we know $\pk{ \sup_{t\ge 0} Y_i(t)>x}=o (\pk{T>x})$ as $x\to\IF$. Therefore, the claim for $\vk M^{(1)}$ is a direct consequence of \nelem{Lem:sum} and \nelem{Lem:compar}. This completes the proof. \QED

%The next lemma is about an auxiliary random vector $\widetilde{\vk U}^{\vk y}=( \widetilde{U_1}^{\vk y}, \widetilde{ U_2}^{\vk y})$ defined  as
%\BQN\label{eq:Uny}
%\widetilde{U_i}^{\vk y}:=\left\{
%  \begin{array}{ll}
% M^{(1)}_i, & \hbox{if \ $M^{(1)}_i> y_1$}, \\
%   U^{(1)}_i, & \hbox{if \ $-y_2<   U^{(1)}_i\le M^{(1)}_i\le y_1$}, \\
%    -y_2, & \hbox{if \ $M^{(1)}_i\le y_1,   U^{(1)}_i\le -y_2 $},
%  \end{array}
%\right.\ \ i=1,2,
%\EQN
%for some $y_1, y_2>0$.
%In the following,  the superscript $\vk y$  will be compressed in $\widetilde{\vk U}^{\vk y}$ for notational simplicity.

%\BEL\label{Lem:UM}
% For any $y_1, y_2>0$, let  $\widetilde{\vk U}=\widetilde{\vk U}^{\vk y}$ be defined in \eqref{eq:Uny}. We have
%\BQN\label{eq:MUUy}
%\pk{\abs{\widetilde{\vk U}} >x} \sim \pk{\abs{\vk M^{(1)}} >x} \sim \pk{\abs{\vk U^{(1)}} >x}\sim  \pk{\abs{\vk{\mathcal{T}}}>x}, \ \ \ \ x\to\IF.
%\EQN
% Furthermore, $\widetilde{\vk U}$ is regularly varying with index $\lambda$ and limiting measure $\mu$.
%\EEL

\prooftheo{Thm:Qu} First,  note that, for any $\vk a, \vk c>\vk 0$,
%\BQNY
%\int_0^\IF \mu( (v\vk c+\vk a, \vk \IF])  dv<\IF.
%\EQNY
%In fact,
by the homogeneous property of $\mu$, %we have
\BQN\label{int-mu-fi}
\int_0^\IF \mu( (v\vk c+\vk a, \vk \IF])  dv\leq \mu( (\vk a, \vk \IF]) +
\int_1^\IF v^{-\lambda}\mu( (\vk c+\vk a/v, \vk \IF])  dv
\leq \mu( (\vk a, \vk \IF]) + \frac{1}{\lambda-1}\mu( (\vk c, \vk \IF]).
\EQN
For simplicity we denote
$\vk W^{(n)}:=\sum_{i=1}^n  \vk U^{(i)}.$  We consider the lower bound, for which we adopt a standard technique of "one big jump" (see \cite{PZ07}).   Informally speaking, we choose an event on which $ \vk W^{(n-1)}+ \vk M^{(n)}, n\ge 1,$ behaves in a typical way up to some time $k$ for which $\vk M^{(k+1)}$ is large.
%%%%%%%%%%
\COM{
For any small $\delta>0$, large $K>0$, and any $n\ge 1$, we define the event
\BQNY
E_n:=E_n(\delta, K)=\{\vk W^{(n-1)}\in (-k(\vk c+\delta\vk 1)-K\vk 1, -k(\vk c-\delta\vk 1)+K\vk 1),\ k\le n \}.
\EQNY
Additionally, consider the event
$$
F_n:=F_n(\delta, K)=\{\vk M^{(k)} <k\delta \vk 1 +K \vk 1,\  k\le n\}.
$$
}%%%%%%%%%%%%%%%%%5
Let $\delta, \vn$ be small positive numbers. By the Weak Law of Large Numbers, we can choose large $K=K_{\vn,\delta}$ so that
\BQNY
\pk{\vk W^{(n)}>-n(1+\vn )\vk c -K\vk 1}>1-\delta,\ \ \ \ n=1,2,\cdots.
\EQNY
For any $u>0$, we have  %(here $\cap_{i=1}^1=\Omega$)
\BQNY
Q(u)&=&\pk{\exists n\ge 1,\  \vk W^{(n-1)} + \vk M^{(n)} >\vk a u}\\
&=&\pk{\vk M^{(1)} >\vk a u}+\sum_{k\ge 1}\pk{\cap_{n=1}^k (\vk W^{(n-1)} + \vk M^{(n)} \not> \vk a u), \vk W^{(k)} + \vk M^{(k+1)} >\vk a u}\\
&\ge &\pk{\vk M^{(1)} >\vk a u}+\sum_{k\ge 1}\mathbb P\Big\{\cap_{n=1}^k (\vk W^{(n-1)} + \vk M^{(n)} \not> \vk a u), \vk W^{(k)}>-k(1+\vn )\vk c-K\vk 1,\\
&& \qquad \ \ \ \ \ \ \ \ \ \ \ \ \ \ \ \ \ \ \ \ \qquad \ \ \ \ \ \ \ \ \qquad \ \ \ \ \ \ \ \ \ \ \ \ \ \ \ \ \ \  \ \ \ \  \ \ \ \ \ \ \ \ \ \ \ \ \vk M^{(k+1)} >\vk a u +k(1+\vn )\vk c+K\vk 1\Big\}\\
&\ge &\pk{\vk M^{(1)} >\vk a u}+\sum_{k\ge 1}\LT(1-\delta-\pk{\cup_{n=1}^k (\vk W^{(n-1)} + \vk M^{(n)} > \vk a u)}\RT) \pk{\vk M^{(k+1)} >\vk a u +k(1+\vn )\vk c+K\vk 1}\\
&\ge &\LT(1-\delta-Q(u)\RT)  \sum_{k\ge 0}\pk{\vk M^{(1)} >\vk a u +k(1+\vn )\vk c+K\vk 1}\\
&\ge &\frac{\LT(1-\delta-Q(u)\RT) }{1+\vn} \int_0^\IF \pk{\vk M^{(1)} >\vk a u + v\vk c+K\vk 1} dv.
\EQNY
For $u$ sufficiently large such that $\vn u>K$, we have
\BQNY
Q(u) \ge \frac{\LT(1-\delta-Q(u)\RT) }{1+\vn} \int_0^\IF \pk{\vk M^{(1)} > (\vk a+\vn \vk 1) u + v\vk c} dv.
\EQNY
Rearranging the above inequality and using a change of variable, we obtain %(more precise proof is needed)
\BQN\label{eq:Qu}
Q(u) \ge  \frac{(1-\delta)  u \int_0^\IF \pk{\vk M^{(1)} >u(\vk a+\vn \vk 1+ v\vk c) } dv }{1+\vn + \int_0^\IF \pk{\vk M^{(1)} >(\vk a+\vn \vk 1) u+ v\vk c} dv}, %\\
%&\sim & \frac{ 1-\delta }{1+\vn }   \int_0^\IF \pk{\vk M^{(1)} >(\vk a+\vn \vk 1) u} dv.
\EQN
and thus by \nelem{lem:MU} and Fatou's lemma
\BQNY
\liminf_{u\to\IF}\frac{Q(u)}{ u \pk{\abs{\wtT }> u}} \ge \frac{1-\delta}{1+\vn}\int_0^\IF \mu( (\vk a +\vn \vk 1+v\vk c, \vk \IF])  dv.
\EQNY
Since $\vn$ and $\delta$ are arbitrary, and by \eqref{int-mu-fi} the integration on the right hand side is finite,  taking $\vn\to0, \delta\to0$ and applying dominated convergence theorem yields
\BQNY\label{Qu:lowbound}
\liminf_{u\to\IF}\frac{Q(u)}{ u \pk{\abs{\wtT }> u}} \ge \int_0^\IF \mu( (\vk a +v\vk c, \vk \IF])  dv.
\EQNY

Next, we consider the asymptotic upper bound. Let $y_1, y_2>0$ be given. We shall  construct an auxiliary random walk $\widetilde{\vk W}^{(n)}, n\ge 0$, with $\widetilde{\vk W}^{(0)}=0$ and %. For $n\ge 1,$ let
\COM{ %%%%%%%%%%5
\BQNY %\label{eq:Uny}
U_i^{(n),y}=\left\{
  \begin{array}{ll}
  U^{(n)}_i, & \hbox{if \ $M^{(n)}_i\le y$;} \\
 M^{(n)}_i, & \hbox{if \ $M^{(n)}_i> y$},
  \end{array}
\right.
\qquad i=1,2
\EQNY
}%%%%%%%%%%%%%%5
 $\widetilde{\vk W}^{(n)}=\sum_{i=1}^n \widetilde{\vk U}^{(i)}, n\ge 1$, where $\widetilde{\vk U}^{(n)} =(\widetilde{U}_1^{(n)}, \widetilde{U}_2^{(n)})$ is given by
\BQNY %\label{eq:Uny}
\widetilde{U}_i^{(n)}=\left\{
  \begin{array}{ll}
 M^{(n)}_i, & \hbox{if \ $M^{(n)}_i> y_1$}; \\
   U^{(n)}_i, & \hbox{if \ $-y_2<   U^{(n)}_i\le M^{(n)}_i\le y_1$}; \\
    -y_2, & \hbox{if \ $M^{(n)}_i\le y_1,   U^{(n)}_i\le -y_2 $},
  \end{array}
\right.
\qquad i=1,2.
\EQNY
 Obviously, $\vk W^{(n)}\le \widetilde{\vk W}^{(n)}$ for any $n\ge 1.$ Furthermore, one can show that % (see also \cite{PZ07})
\BQNY
M_i^{(n)} %&\le & \widetilde{U}_i^{(n)} +(y_1-\widetilde{U}_i^{(n)}) I_{(-y_2<   U^{(n)}_i\le M^{(n)}_i\le y_1)} + (y_1-\widetilde{U}_i^{(n)}) I_{(M^{(n)}_i\le y_1,   U^{(n)}_i\le -y_2) }\\
&\le & \widetilde{U}_i^{(n)} + (y_1+y_2).
\EQNY
%\xP{\{The first inequality is not right! But the last inequality is true.\}}
Then,
\BQNY
\vk W^{(n-1)} + \vk M^{(n)}  \le \widetilde{\vk W}^{(n)} + (y_1+y_2) \vk 1,\ \ \ n\ge 1.
\EQNY
Thus, for any $\vn>0$  and sufficiently large $u$,
\BQNY
% \pk{\exists n\ge 1,\ \vk W^{(n-1)} + \vk M^{(n)} >\vk a u}
Q(u)&\le& \pk{\exists n\ge 1,\ \widetilde{\vk W}^{(n)}   >\vk a u - (y_1+y_2) \vk 1}\\
&\le& \pk{\exists n\ge 1,\ \widetilde{\vk W}^{(n)}   >(\vk a -\vn \vk 1) u }.
\EQNY
Define $c_{y_1,y_2}=-\E{\widetilde{\vk U}^{(1)}}$. Since $\lim_{y_1, y_2\to\IF} \vk c_{y_1,y_2}=\vk c$, we have   that for any  $y_1,y_2$ large enough  $\vk c_{y_1,y_2} >\vk 0 $.
It follows from \nelem{lem:MU} and \nelem{Lem:compar} that for any $y_1, y_2>0$, $\widetilde{\vk U}^{(1)}$ is regularly varying with index $\lambda$ and limiting measure $\mu$, and $\pk{ \abs{\widetilde{\vk U}^{(1)} }>u}\sim\pk{\abs{\wtT}>u}$ as $u\to\infty$.
 Then, applying Theorem 3.1 and Remark 3.2 of %Hult, Lindskog, Mikosch and Samorodnitsky (2006)
 \cite{Hult06} we obtain that
\BQNY
\pk{\exists n\ge 1,\  \widetilde{\vk W}^{(n-1)} >(\vk a -\vn \vk 1) u}
&\sim& u\pk{ \abs{\widetilde{\vk U}^{(1)} }>u} \int_0^\IF \mu ((\vk c_{y_1,y_2} v+\vk a-\vn \vk 1, \vk\IF]) dv\\
&\sim&  u     \pk{\abs{\wtT}>u} \int_0^\IF \mu ((\vk c_{y_1,y_2} v+\vk a-\vn \vk 1, \vk\IF]) dv.
%&\sim & \int_0^\IF \pk{\vk M^{(1)} >\vk a u + v\vk c^y } dv,
\EQNY
Consequently, the claimed asymptotic upper bound is obtained by letting $\vn\to0$, $y_1, y_2\to\IF$.  The proof is complete. \QED

\section{Appendix }

This section includes some results on the regularly varying random vectors.% followed by the technical proofs of Lemma \ref{lem:MU}  and Lemma \ref{Lem:UM}.

\BEL\label{Lem:1}
 Let $\vk{\mathcal{T}}>\vk 0$ be a regularly varying random vector with index $\alpha$ and limiting measure $\nu$, and let $x_i(u), 1\le i\le n$ be increasing  (to infinity)  functions     such that for some $1\le m \le n,$ $x_1(u) \sim \cdots \sim x_m(u)$, and $x_j(u)=o(x_{1}(u))$ for all $j=m+1,\cdots, n$. Then, for any $\vk a>\vk 0$, %$a_i>0, 1\le i\le n$,
\BQNY
\pk{\cap_{i=1}^n  (\mathcal{T}_i > a_i x_i(u)) }\sim \pk{\cap_{i=1}^m  (\mathcal{T}_i > a_i x_1(u)) } \sim \nu([\vk a_{m,0}, \vk\IF])\ \pk{\abs{\vk{\mathcal{T}}} > x_1(u)}
\EQNY
holds as $u\to \IF$, with $\vk a_{m,0}=(a_1,\cdots, a_m,  0, \cdots, 0)$.
\EEL

\prooflem{Lem:1} Obviously, for any small enough $\vn>0$ we have that when  $u$ is sufficiently large
\BQNY
\pk{\cap_{i=1}^n  (\mathcal{T}_i > a_i x_i(u)) }&\le& \pk{\cap_{i=1}^m  (\mathcal{T}_i > (a_i-\vn) x_1(u)), \cap_{i=m+1}^n  (\mathcal{T}_i  > 0  )} \\
&\sim &\nu([\vk a_{-\vn}, \vk\IF])\ \pk{\abs{\vk{\mathcal{T}}} > x_1(u)}, %\pk{\cap_{i=1}^m  (\mathcal{T}_i > (a_i+\vn)  x_1(u)) }
\EQNY
where  $\vk a_{-\vn}=(a_1-\vn,\cdots, a_m-\vn,  0, \cdots, 0)$.
Next, for any small enough $\vn>0$ we have that when  $u$ is sufficiently large
\BQNY
\pk{\cap_{i=1}^n  (\mathcal{T}_i > a_i x_i(u)) }  & \ge & \pk{\cap_{i=1}^m  (\mathcal{T}_i > (a_i+\vn) x_1(u)), \cap_{i=m+1}^n  (\mathcal{T}_i  >a_i (\vn x_1(u)))}\\
&\sim&\nu([\vk a_{\vn+}, \vk\IF]) \pk{\abs{\vk{\mathcal{T}}} > x_1(u)}
\EQNY
with $\vk a_{\vn+}=(a_1+\vn,\cdots, a_m+\vn,  a_{m+1} \vn, \cdots, a_{n}\vn)$.
\COM{\BQNY
&&\lim_{u\to\IF}\frac{\pk{\cap_{i=1}^n  (\mathcal{T}_i > a_i u^{\lambda_i}) }  }{\pk{\abs{\vk{\mathcal{T}}} > u^{\lambda_1}} }   \ge\lim_{u\to\IF}\frac{\pk{\cap_{i=1}^m  (\mathcal{T}_i > a_i u^{\lambda_1}), \cap_{i=m+1}^n  (T_2 >a_i (\vn u^{\lambda_1}))} }{\pk{\abs{\vk{\mathcal{T}}} > u^{\lambda_1}}} \\
&&=\nu([\vk a_{\vn}, \vk\IF]),
\EQNY}
Letting $\vn \to 0$, the claim follows by the continuity of  $\nu([\vk a_{\vn\pm}, \vk\IF])$ in $\vn$. The proof is complete. \QED

\BEL\label{Lem:2}
 Let $\vk{\mathcal{T}}$, $a_i$'s, $x_i(u)'s$ and $\vk a_{m,0}$ be the same as in \nelem{Lem:1}. Further, consider ${\vk \eta}=(\eta_1, \cdots, \eta_n)$ to be  an independent of $\vk{\mathcal{T}}$ nonnegative random vector %with arbitrarily dependent components satisfying
 such that $\max_{1\le i\le n}\E{\eta_i^{\alpha+\delta}}<\IF$ for some $\delta>0$.  Then, %for any $\vk a>\vk 0,$ %$a_i>0, 1\le i\le n$
\BQNY
\pk{\cap_{i=1}^n  (\mathcal{T}_i \eta_i> a_i x_i(u)) }\sim \pk{\cap_{i=1}^m  (\mathcal{T}_i \eta_i> a_i x_1(u)) }\sim \widehat\nu([\vk a_{m,0}, \vk\IF])\ \pk{\abs{\vk{\mathcal{T}}} > x_1(u)}
\EQNY
holds as $u\to \IF$, where $\widehat{\nu}(K)=\E{\nu(\vk \eta^{-1} K)}$, with $\vk \eta^{-1} K=\{(\eta_1^{-1}b_1,\cdots, \eta_n^{-1}b_n), (b_1,\cdots,b_n)\in K\}$ for any $K\subset  \mathcal{B}([0,\IF]^n \setminus \{\vk 0\})$.
\EEL

\prooflem{Lem:2}
It follows directly from Lemma 4.6 of \cite{Jessen06}   (see also Proposition A.1 of \cite{Basrak02}) that
%By Lemma 3.1 of Konstantinides and Li (2016), we conclude that
%$  \vk{\mathcal{T}} \vk \xi :=(\mathcal{T}_1\xi_1, \cdots, \mathcal{T}_n\xi_n)$ is regularly varying with index $\alpha$, and that
the second asymptotic equivalence holds.
The first claim follows from the same arguments as in \nelem{Lem:1}. \QED

\BEL \label{Lem:sum}
Assume $\vk X\in \R^n$ is regularly varying with index $\alpha$ and limiting measure $\mu$, $\vk A$ is a random
$n\times d$ matrix independent of random vector
$\vk Y\in \R^d$. If $0<\E{\norm{\vk A}^{\alpha+\delta}}<\infty$ for some $\delta>0$, with $\norm{\cdot}$ some matrix norm % (e.g., maximum norm)
and
\BQN\label{YnegX}
\pk{\abs{\vk Y}>x} = o\LT(\pk{\abs{\vk X}>x}\RT),\ \ \ \ x\to\IF,
\EQN
then, $\vk X+\vk A \vk Y$ is regularly varying with index $\alpha$ and limiting measure $\mu$, and
\BQNY
\pk{\abs{\vk X+\vk A \vk Y}>x} \sim \pk{\abs{\vk X}>x},\ \ \ \ x\to\IF.
\EQNY
\EEL

\prooflem{Lem:sum} By Lemma 3.12 of \cite{Jessen06}, it suffices to show that
\BQN\label{AYnelX}
\pk{\abs{\vk A \vk Y}>x} = o\LT(\pk{\abs{\vk X}>x}\RT),\ \ \ \ x\to\IF.
\EQN
 Defining $g(x)=x^{\frac{\alpha+\delta/2}{\alpha+\delta}}, x\ge 0$, we have
\BQN\label{eq:xiS}
\pk{\abs{\vk A \vk Y}>x} \le \pk{\norm{\vk A}\abs{\vk Y}>x} \le
 \int_0^{g(x)}  \pk{\abs{\vk Y}>x/t} \pk{\norm{\vk A} \in dt} +\pk{\norm{\vk A} >g(x)}.
\EQN
Due to \eqref{YnegX}, for arbitrary $\vn>0$,
\BQNY
\int_0^{g(x)}\pk{\abs{\vk Y}>x/t}\pk{\norm{\vk A} \in dt}
\le \vn \int_0^{g(x)} \pk{\abs{\vk X}>x/t}\pk{\norm{\vk A} \in dt},
\EQNY
hold for large enough $x$. Furthermore, by Potter's Theorem (see, e.g.,  Theorem 1.5.6 of \cite{bingham1989regular}), we have
\BQNY
\frac{\pk{\abs{\vk X}>x/t}}{\pk{\abs{\vk X}>x}}\le I_{(t\le 1)} + 2 t^{\alpha+\delta}I_{(1<t\le g(x))}, \quad t\in (0, g(x))
\EQNY
holds for sufficiently large $x$, and thus by the dominated convergence theorem,
\BQN\label{eq:gx0}
\lim_{x\to\IF}   \int_0^{g(x)}   \frac{\pk{\abs{\vk Y}>x/t}}{ \pk{\abs{\vk X}>x}}  \pk{\norm{\vk A} \in dt}
\le \lim_{x\to\IF}   \int_0^{g(x)}   \frac{\vn\pk{\abs{\vk X}>x/t}}{ \pk{\abs{\vk X}>x}}  \pk{\norm{\vk A} \in dt}
= \vn \E{\norm{\vk A}^{\alpha}}.
\EQN
%for any $\vn>0$. Letting $\vn\to0$, we obtain
%\BQN\label{eq:gx0}
%\lim_{x\to\IF}   \int_0^{g(x)}   \frac{\pk{S>x/u}}{ \pk{T>x}}  \pk{\xi \in du}  =0.
%\EQN
Moreover, Markov inequality implies that
\BQN\label{eq:gx1}
 \lim_{x\to\IF}  \frac{\pk{\norm{\vk A}>g(x)}}{ \pk{\abs{\vk X}>x}} \le  \lim_{x\to\IF}  \frac{ \E{\norm{\vk A}^{\alpha+\delta}}}{ g(x)^{\alpha+\delta}\pk{\abs{\vk X}>x}}=0.
\EQN
Therefore, the claim \eqref{AYnelX} follows from \eqref{eq:xiS}-\eqref{eq:gx1} and the arbitrariness of $\vn$. This completes the proof. \QED

\BEL \label{Lem:compar}
Assume $\vk X , \vk Y \in \R^n$ are regularly varying with same index $\alpha$ and same limiting measure $\mu$. Moreover, if $\vk X \geq \vk Y$ and
$\pk{\abs{\vk X}>x} \sim \pk{\abs{\vk Y}>x}$ as  $x\to\IF$,
%\BQN\label{eq:teqixy}
%\pk{\abs{\vk X}>x} \sim \pk{\abs{\vk Y}>x}\quad \textrm{as}\  x\to\IF,
%\EQN
 then for any random vector $\vk Z$ satisfying  $\vk X \geq \vk Z \geq \vk Y$, $\vk Z$ is regularly varying with index $\alpha$ and limiting measure $\mu$, and $\pk{\abs{\vk Z}>x} \sim \pk{\abs{\vk X}>x}$ as $x\to\IF$.
\EEL

\prooflem{Lem:compar} We only prove the claim for $n=2$, a similar argument can be used to verify the claim for $n\geq3$.
 For any $x>0$, define a measure $\mu_x$ as
\BQNY\label{def:mux}
\mu_x(A)=:\frac{\pk{x^{-1}{\vk Z}\in A}}{\pk{\abs{\vk X}>x}},\quad A\in \mcB(\barR_0^2).
\EQNY
We shall  show that
\BQN \label{eq:muxmu}
\mu_x \stackrel{v}{\longrightarrow} \mu, \quad  x\to\infty.
\EQN
Given that %\eqref{eq:muxmu}
the above
 is established, by letting $A=\{\vk x: |\vk x|>1\}$ (which is relatively compact and satisfies $\mu(\partial A)=0$),  we have  $\mu_x(A)\to\mu(A)=1$ as $x\to\infty$ and thus $\pk{|\vk Z|>x}\sim\pk{|\vk X|>x}$. Furthermore, by substituting the denominator in the definition of $\mu_x$ by $\pk{|\vk Z|>x}$, we conclude that
 \BQNY
\frac{\pk{x^{-1}{\vk Z}\in \cdot}}{\pk{|\vk Z|>x}}\stackrel{v}{\longrightarrow} \mu(\cdot), \quad  x\to\infty,
\EQNY
showing that $\vk Z$ is regularly varying with index $\alpha$ and limiting measure $\mu$.

%We first show that $\{\mu_x\}_{x>0}$ is relatively compact.
% It suffices to show that
%\BQN\label{eq:Mab}
%\lim_{x\to\IF}\frac{\pk{x^{-1} \vk Z\in (\vk a, \vk b]  }}{\pk{\abs{\vk Z}>x}}= \mu((\vk a, \vk b])
%\EQN
%for all rectangle sets $(\vk a, \vk b] \in \barR_0^2$.
Now it remains to prove \eqref{eq:muxmu}. To this end, we define a set $\mathcal D$ consisting of all sets in $ \barR_0^2$ that are of the following form:
\BQNY
&&a):\quad (a_1,\IF]\times [a_2,\IF],\ \ \ \ a_1>0, a_2\in \R,\\
&&b):\quad [-\IF, a_1]\times (a_2,\IF],\ \ \ \ a_1\in \R, a_2>0,\\
&&c):\quad [-\IF, a_1)\times [-\IF, a_2],\ \ \ \ a_1<0, a_2\in \R,\\
&&d):\quad [a_1,\IF]\times [-\IF, a_2),\ \ \ \ a_1\in \R, a_2<0.
\EQNY
%One can easily check that any rectangle $(\vk a, \vk b] \in \barR_0^2$ can be obtained from a finite number of  sets in $\mathcal D$ by operating union, intersection, difference or complementary. Thus, in order to derive \eqref{eq:Mab} it is sufficient to show that
Note that every $A\in\mathcal D$ is relatively compact and satisfies $\mu(\partial A)=0$.
We first show that
\BQN\label{eq:mux-mu}
%\lim_{x\to\IF}\frac{\pk{x^{-1} \vk Z\in K  }}{\pk{\abs{\vk Z}>x}}= \mu( K)
\lim_{x\to\IF} \mu_x(A)= \mu( A),\ \ \ \ \forall A\in \mathcal D.
\EQN
%holds for all sets $A\in \mathcal D$.
 If $A=(a_1,\IF]\times (a_2,\IF]$ or $A=(a_1,\IF]\times [a_2,\IF]$ with $a_i\in \R$ and at least one $a_i>0, i=1,2$, or $A=\barR\times (a_2,\IF]$ with some $a_2>0$, by the order relations of $\vk X,\vk Y,\vk Z$, we have for any $x>0$
\BQN\label{ineq:mux}
\frac{\pk{x^{-1} \vk Y\in A  }}{\pk{\abs{\vk X}>x}}\le \mu_x(A) \le \frac{\pk{x^{-1} \vk X\in A  }}{\pk{\abs{\vk X}>x}}.
\EQN
Letting $x\to\infty$, using the regularity properties as supposed for $\vk X$ and $\vk Y$, and then appealing to Proposition 3.12(ii) in \cite{Res87}, we verify  \eqref{eq:mux-mu} for case a). If $A=[-\IF, a_1]\times (a_2,\IF]$ with some $a_1\in \R, a_2>0$, then we have
\BQNY
\mu_x( A)=\mu_x(\barR\times(a_2,\IF]) - \mu_x((a_1,\IF]\times (a_2,\IF]),
\EQNY
and thus by the convergence in case a),
\BQNY
\lim_{x\to\IF}\mu_x(A)= \mu( \barR\times  (a_2,\IF])- \mu((a_1,\IF]\times  (a_2,\IF]) =\mu(A),
\EQNY
this validates  \eqref{eq:mux-mu} for case b).
%\BQNY
% \pk{x^{-1}  \vk Z\in A  } = \pk{x^{-1} \vk Z\in  \barR\times  (a_2,\IF]} - \pk{x^{-1} \vk Z\in   (a_1,\IF]\times  (a_2,\IF]}.
%\EQNY
%Similarly as the previous case we see that  \eqref{eq:Mab2}  is valid for sets of the form $(a_1,\IF]\times (a_2,\IF],$ for any  $a_1\in \barR, a_2>0$. Thus,
%\BQNY
%\lim_{x\to\IF}\frac{\pk{x^{-1} \widetilde{\vk U}\in K  }}{\pk{\abs{\widetilde{\vk U}}>x}}= \mu( \barR\times  (a_2,\IF])- \mu((a_1,\IF]\times  (a_2,\IF]) =\mu(K).
%\EQNY
If $A= [-\IF, a_1)\times [-\IF, a_2]$ or $A= [-\IF, a_1)\times [-\IF, a_2)$  with $a_i\in \R$ and at least one $a_i<0,  i=1,2$, or $A= \barR\times [-\IF, a_2)$ with some $a_2<0$, then we get a similar formula as \eqref{ineq:mux} with the reverse inequalities.
% \BQNY
%\frac{\pk{x^{-1} \vk Y\in A  }}{\pk{\abs{\vk X}>x}}\ge \mu_x(A) \ge \frac{\pk{x^{-1} \vk X\in A  }}{\pk{\abs{\vk X}>x}}.
%\EQNY
%\BQNY
%0\le \pk{x^{-1} \widetilde{\vk U}\in K  } \le \pk{x^{-1} \vk U^{(1)}\in K  }.
%\EQNY
%The claim in \eqref{eq:Mab2} follows since now $\mu(K)=0$. Lastly,
If $A=[a_1,\IF]\times [-\IF, a_2)$ with some $a_1\in \R, a_2<0,$ then
\BQNY
\mu_x(A)=\mu_x(\barR\times  [-\IF, a_2))-\mu_x([-\IF, a_1)\times  [-\IF, a_2)).
\EQNY
Therefore, similarly as the proof for the cases a)-b), one can establish \eqref{eq:mux-mu} for the   cases c) and d).
%  we see that \eqref{eq:Mab2} holds since  $\mu(K)=0$ for such $K$.

Next, let $f$ defined on $ \barR_0^2$ be any positive, continuous function with compact support. We see that the support of $f$ is contained in $[\vk a,\vk b]^c$ for some $\vk a<\vk 0<\vk b$.
%Since we can obtain $(\vk a,\vk b]^c$ by finite set operation on $\mathcal D$ and on which have convergence \eqref{eq:mux-mu}, therefore,
Note that
\BQNY
[\vk a,\vk b]^c=(b_1,\infty]\times[a_2,\infty]\cup[-\infty,b_1]\times(b_2,\infty]\cup
[-\infty,a_1)\times[-\infty,b_2]\cup[a_1,\infty]\times[-\infty,a_2)
=:\bigcup_{i=1}^4A_i,
\EQNY
where $A_i$'s are sets of the form  a)-d) respectively, and thus   \eqref{eq:mux-mu} holds for these $A_i$'s. Therefore,
\BQNY
\sup_{x>0}\mu_x(f)\leq \sup_{\vk z\in\barR_0^2}f(\vk z) \cdot \sup_{x>0}\mu_x([\vk a,\vk b]^c)\leq \sup_{\vk z\in\barR_0^2}f(\vk z) \cdot \sum_{i=1}^4\sup_{x>0}\mu_x(A_i)<\infty,
\EQNY
which by Proposition 3.16 of \cite{Res87} implies that $\{\mu_x\}_{x>0}$ is a vaguely relatively compact subset of the metric space consisting of %. Moreover, by Proposition 3.17 of \cite{Res87}, all
all the nonnegative Radon measures on $(\barR_0^2, \mcB(\barR_0^2))$. %can be metrizable as a complete, separable  metric space.
If $\mu_0$ and $\mu_0'$ are two subsequential vague limits of $\{\mu_x\}_{x>0}$ as $x\to\infty$, then by \eqref{eq:mux-mu} % and Proposition 3.12(ii) in \cite{Res87},
we have $\mu_0(A)=\mu_0'(A)$ for any $A\in \mathcal D$. Since any rectangle   in $\barR_0^2$ can be obtained from a finite number of  sets in $\mathcal D$ by operating union, intersection, difference or complementary, and these rectangles constitutes a $\pi$-system and generate the $\sigma$-field $\mcB(\barR_0^2)$, we get $\mu_0=\mu_0'$ on $ \barR_0^2 $. Consequently, \eqref{eq:muxmu} is valid, and thus the proof is complete.
% $\mu_x$ converges vaguely to $\mu$ as $x\to\infty$.
%\BQNY
%\mu_x \stackrel{v}{\longrightarrow} \mu, \quad  x\to\infty.
%\EQNY
 % This completes the proof.
 \QED

\bigskip
{\bf Acknowledgement}:
%We are thankful to the editor and the referee for their constructive suggestions which have led to a significant improvement of the manuscript.
The research of Xiaofan Peng is partially supported by  National Natural Science Foundation of China (11701070,71871046).

\bibliographystyle{ieeetr}

 \bibliography{gausbibR}

\newcommand{\nosort}[1]{}\def\polhk#1{\setbox0=\hbox{#1}{\ooalign{\hidewidth
  \lower1.5ex\hbox{`}\hidewidth\crcr\unhbox0}}}
  \def\polhk#1{\setbox0=\hbox{#1}{\ooalign{\hidewidth
  \lower1.5ex\hbox{`}\hidewidth\crcr\unhbox0}}} \def\cprime{$'$}
  \def\cprime{$'$}
\begin{thebibliography}{10}

\bibitem{debicki2004supremum}
K.~D\c{e}bicki, B.~Zwart, and S.~Borst, ``The supremum of a {G}aussian process
  over a random interval,'' {\em Statistics and Probability letters}, vol.~68,
  no.~3, pp.~221--234, 2004.

\bibitem{arendarczyk2011asymptotics}
M.~Arendarczyk and K.~D\c{e}bicki, ``Asymptotics of supremum distribution of a
  {G}aussian process over a {W}eibullian time,'' {\em Bernoulli}, vol.~17,
  no.~1, pp.~194--210, 2011.

\bibitem{arendarczyk2011exact}
M.~Arendarczyk and K.~D\c{e}bicki, ``Exact asymptotics of supremum of a
  stationary {G}aussian process over a random interval,'' {\em Statistics and
  Probability Letters}, 2011.

\bibitem{TH13}
Z.~Tan and E.~Hashorva, ``Exact asymptotics and limit theorems for supremum of
  stationary chi-processes over a random interval,'' {\em Stochastic Processes
  and Their Applications}, vol.~123, pp.~2983--2998, 2013.

\bibitem{DEJ14}
K.~D{\polhk{e}}bicki, E.~Hashorva, and L.~Ji, ``Tail asymptotics of supremum of
  certain {G}aussian processes over threshold dependent random intervals,''
  {\em Extremes}, vol.~17, no.~3, pp.~411--429, 2014.

\bibitem{Are17}
M.~Arendarczyk, ``On the asymptotics of supremum distribution for some
  interated processes,'' {\em Extremes}, vol.~20, pp.~451--474, 2017.

\bibitem{DP20}
K.~D{\polhk{e}}bicki and X.~Peng, ``Sojourns of stationary {G}aussian processes
  over a random interval,'' {\em Preprint,
  https://arxiv.org/pdf/2004.12290.pdf}, 2020.

\bibitem{Debicki10}
K.~D{\polhk{e}}bicki, K.~Kosi{\'n}ski, M.~Mandjes, and T.~Rolski, ``Extremes of
  multi-dimensional {G}aussian processes,'' {\em Stochastic Processes and their
  Applications}, vol.~120, no.~12, pp.~2289--2301, 2010.

\bibitem{DHJT15}
K.~D\c{e}bicki, E.~Hashorva, L.~Ji, and K.~Tabi\'{s}, ``Extremes of
  multi-dimensional {G}aussian processes: {E}xact asymptotics,'' {\em
  Stochastic Processes and Their Applications}, vol.~125, pp.~4039--4065, 2015.

\bibitem{AP19}
J.-M. Azais and V.~Pham, ``Geometry of conjunction set of smooth stationary
  {G}aussian fields,'' {\em Preprint, https://arxiv.org/abs/1909.07090v1},
  2019.

\bibitem{Pham20}
V.~Pham, ``Conjunction probability of smooth centered {G}aussian processes,''
  {\em Acta Math Vietnam, https://doi.org/10.1007/s40306-019-00351-4}, 2020.

\bibitem{DHW20}
K.~D{\polhk{e}}bicki, E.~Harshorva, and L.~Wang, ``Extremes of
  multi-dimensional {G}aussian processes,'' {\em Stochastic Processes and their
  Applications}, 2020.

\bibitem{WF00}
K.~Worsley and K.~Friston, ``A test for a conjunction,'' {\em Statistics and
  Probability Letters}, vol.~47, no.~2, pp.~135--140, 2000.

\bibitem{MR2145669}
V.~Piterbarg and B.~Stamatovich, ``Rough asymptotics of the probability of
  simultaneous high extrema of two {G}aussian processes: the dual action
  functional,'' {\em Uspekhi Mat. Nauk}, vol.~60, no.~1(361), pp.~171--172,
  2005.

\bibitem{DJT20}
K.~D{\polhk{e}}bicki, L.~Ji, and T.~Rolski, ``Exact asymptotics of
  component-wise extrema of two-dimensional {B}rownian motion,'' {\em Extremes,
  https://doi.org/10.1007/s10687-020-00387-y}, 2020.

\bibitem{DHK20}
K.~D{\polhk{e}}bicki, E.~Harshorva, and K.~Krystecki, ``Finite-time ruin
  probability for correlated {B}rownian motions,'' {\em Preprint,
  https://arxiv.org/pdf/2004.14015.pdf}, 2020.

\bibitem{HH19}
R.~van~der Hofstad and H.~Honnappa, ``Large deviations of bivariate {G}aussian
  extrema,'' {\em Queueing Systems}, vol.~93, pp.~333--349, 2019.

\bibitem{AsmAlb10}
S.~Asmussen and H.~Albrecher, {\em Ruin probabilities}.
\newblock Advanced Series on Statistical Science \& Applied Probability, 14,
  World Scientific Publishing Co. Pte. Ltd., Hackensack, NJ, second~ed., 2010.

\bibitem{JR18}
L.~Ji and S.~Robert, ``Ruin problem of a two-dimensional fractional {B}rownian
  motion risk process,'' {\em Stochastic Models}, vol.~34, no.~1, pp.~73--97,
  2018.

\bibitem{BPS01}
O.~E. Barndorff-Nielsen, J.~Pedersen, and K.-I. Sato, ``Multivariate
  subordination, self-decomposability and statility,'' {\em Advances in Applied
  Probability}, vol.~33, pp.~160--187, 2001.

\bibitem{LS10}
E.~Luciano and P.~Semeraro, ``Multivariate time changes for {L}\'evy asset
  models: {C}haracterization and calibration,'' {\em Journal of Computational
  and Applied Mathematics}, vol.~233, pp.~1937--1953, 2010.

\bibitem{Kim12}
Y.~S. Kim, ``The fractional multivariate normal tempered stable process,'' {\em
  Applied Mathematics Letters}, vol.~25, pp.~2396--2401, 2012.

\bibitem{HKR98}
H.~He, W.~P. Keirstead, and J.~Rebholz, ``Double lookbacks,'' {\em Mathematical
  Finance}, vol.~8, no.~3, pp.~201--228, 1998.

\bibitem{PZ07}
Z.~Palmowski and B.~Zwart, ``Tail asymptotics of the supremum of a regenerative
  process,'' {\em Journal of Applied Probability}, vol.~44, no.~2,
  pp.~349--365, 2007.

\bibitem{ZBD05}
B.~Zwart, S.~Borst, and K.~D{\polhk{e}}bicki, ``Subexponential asymptotics of
  hybrid fluid and ruin models,'' {\em The Annals of Applied Probability},
  vol.~15, no.~1A, pp.~500--517, 2005.

\bibitem{KellaWhitt}
O.~Kella and W.~Whitt, ``A storage model with a two-state random evironment,''
  {\em Operations Research}, vol.~40, pp.~257--262, 1992.

\bibitem{Corina20}
C.~Constantinescu, G.~Delsing, M.~Mandjes, and L.~Rojas~Nandayapa, ``A ruin
  model with a resampled environment,'' {\em Scandinavian Actuarial Journal},
  vol.~4, pp.~323--341, 2020.

\bibitem{Rat20}
N.~Ratanov, ``Kac-levy process,'' {\em Journal of Theoretical Probability},
  vol.~33, pp.~239--267, 2020.

\bibitem{Jessen06}
A.~H. Jessen and T.~Mikosch, ``Regular varying functions,'' {\em Publications
  de L'Institut Mathematique}, vol.~80, no.~94, pp.~171--192, 2006.

\bibitem{Res07}
S.~I. Resnick, {\em Heavy-{T}ail {P}henomena. {P}robabilistic and {S}tatistical
  {M}odeling}.
\newblock London: Springer, 2007.

\bibitem{Hult06}
H.~Hult, F.~Lindskog, T.~Mikosch, and G.~Samorodnitsky, ``Functional large
  deviations for multivariate regularly varying random walks,'' {\em The Annals
  of Applied Probability}, vol.~15, no.~4, pp.~2651--2680, 2006.

\bibitem{dieker2005extremes}
A.~Dieker, ``Extremes of {G}aussian processes over an infinite horizon,'' {\em
  Stochastic Processes and their Applications}, vol.~115, no.~2, pp.~207--248,
  2005.

\bibitem{Pit96}
V.~Piterbarg, {\em Asymptotic methods in the theory of {G}aussian processes and
  fields}, vol.~148 of {\em Translations of Mathematical Monographs}.
\newblock Providence, RI: American Mathematical Society, 1996.

\bibitem{ST94}
G.~Samorodnitsky and M.~S. Taqq, ``Stable {N}on-{G}aussian {R}andom
  {P}rocesses,'' {\em The Annals of Applied Probability}, vol.~15, no.~4,
  pp.~2651--2680, 2006.

\bibitem{bingham1989regular}
N.~Bingham, C.~Goldie, and J.~Teugels, {\em Regular variation}, vol.~27.
\newblock Cambridge university press, 1989.

\bibitem{ZN03}
N.~L. Zadi and D.~Nualart, ``Smoothness of the law of the supremum of the
  fractional {B}rownian motion,'' {\em Electronic Communications in
  Probability}, vol.~8, pp.~102--111, 2003.

\bibitem{Mik99}
T.~Mikosch, {\em Regular variation, subexponentiatility and their applications
  in probability theory}.
\newblock in: Lecture Notes for the Workshop "Heavy Tails and Queues", EURANDOM
  Institute, Eindhoven, The Netherlands, 1999.

\bibitem{Basrak02}
B.~Basrak, R.~A. Davis, and T.~Mikosch, ``Regular variation of {GARCH}
  processes,'' {\em Stochastic Processes and their Applications}, vol.~99,
  pp.~95--116, 2002.

\bibitem{Res87}
S.~I. Resnick, {\em Extreme Values, Regular variation, and Point Processes}.
\newblock London: Springer, 1987.

\end{thebibliography}
\end{document}